# CONSISTENT FAMILIES OF BROWNIAN MOTIONS AND STOCHASTIC FLOWS OF KERNELS

BY CHRIS HOWITT AND JON WARREN

*University of Warwick*

Consider the following mechanism for the random evolution of a distribution of mass on the integer lattice **Z**. At unit rate, independently for each site, the mass at the site is split into two parts by choosing a random proportion distributed according to some specified probability measure on $[0,1]$ and dividing the mass in that proportion. One part then moves to each of the two adjacent sites. This paper considers a continuous analogue of this evolution, which may be described by means of a stochastic flow of kernels, the theory of which was developed by Le Jan and Raimond. One of their results is that such a flow is characterized by specifying its $N$ point motions, which form a consistent family of Brownian motions. This means for each dimension $N$ we have a diffusion in $\mathbf{R}^N$, whose $N$ coordinates are all Brownian motions. Any $M$ coordinates taken from the $N$-dimensional process are distributed as the $M$-dimensional process in the family. Moreover, in this setting, the only interactions between coordinates are local: when coordinates differ in value they evolve independently of each other. In this paper we explain how such multidimensional diffusions may be constructed and characterized via martingale problems.

**1. Introduction.** We may say that a pair of Brownian motions $(X_1(t); t \geq 0)$ and $(X_2(t); t \geq 0)$ defined on a common probability space are $\theta$-coupled, where $\theta$ is a positive real parameter if $X_1$ and $X_2$ are both Brownian motions relative to some common filtration and

$$(1.1) \qquad \langle X_1, X_2 \rangle(t) = \int_0^t \mathbf{1}_{(X_1(s) = X_2(s))} \, ds, \qquad t \geq 0,$$

$$(1.2) \qquad L_t^0(X_1 - X_2) = 2\theta \int_0^t \mathbf{1}_{(X_1(s) = X_2(s))} \, ds, \qquad t \geq 0.$$

Here $L_t^0(X_1 - X_2)$ denotes the semimartingale local time accrued by $X_1 - X_2$ at zero by time $t$. A pair of $\theta$-coupled Brownian motions evolve independently of each other when apart, but when they meet there is some interaction, often described as stickiness, which results in their momentarily moving together. The difference of such a pair $(X_1(t) - X_2(t); t \geq 0)$ is a diffusion on $\mathbf{R}$, known as sticky Brownian motion, which is in natural scale and with a speed measure equal to the sum of Lebesgue measure and an atom of size $1/\theta$ at zero. For some previous work on sticky Brownian motion (see Amir [1], Harrison and Lemoine [5] and Warren [17]).

The problem considered in this paper is of describing $\mathbf{R}^N$-valued diffusions $(X(t); t \geq 0)$ with the property that each pair of coordinates $(X_i(t), X_j(t); t \geq 0)$ is a $\theta$-coupled pair of Brownian motions. In fact we are interested in obtaining families of such processes having a natural consistency property. A family consists of a diffusion for each choice of the dimension $N$, and we require that any $M$ coordinates taken from the $N$-dimensional process are distributed as the $M$-dimensional process in the family. We are motivated by work of Le Jan and Raimond [10], who have obtained very general results associating with a consistent family of diffusions a corresponding stochastic flow of kernels. In fact, Le Jan and Raimond [11] as well as Le Jan and Lemaire [9] have already described a consistent family of Brownian motions, of the type we have in mind, by means of Dirichlet forms. Here we will use a martingale problem formulation which will allow us to exhibit many more such families.

Consider the case of $N = 3$. Then, by projecting $(X_1(t), X_2(t), X_3(t); t \geq 0)$ onto the plane $x_1 + x_2 + x_3 = 0$, we obtain a diffusion in $\mathbf{R}^2$ which behaves as standard Brownian motion away from the set of six rays $\{(r\cos\theta, r\sin\theta) \in \mathbf{R}^2 : \theta = 0, \pi/3, 2\pi/3, \pi, 4\pi/3, 5\pi/3\}$. The rays are sticky in a way which is easily described. But the behavior of the process at the origin is more subtle. Ikeda and Watanabe [7] made a careful study of diffusions of this type. The origin is a regular, recurrent point and the time spent there has positive Lebesgue measure. They proved that each excursion made by the process from the origin may be classified as starting along one of the rays and, as a consequence, the behavior of the process at the origin may be characterized by means of a vector of six nonnegative parameters, each of which governs the rate of excursions starting along a particular ray.

The martingale problems which we use to describe our diffusions contain parameters $\theta(k:l)$ where $k$ and $l$ each range through the set of positive integers. The $N$-dimensional process $X$ spends time having a positive Lebesgue measure in certain lower-dimensional subsets of $\mathbf{R}^N$, which we call cells, each corresponding to some collection of equalities and inequalities holding between the coordinates of $X$. The parameter $\theta(k:l)$ has a interpretation similar to that of the parameters in Ikeda and Watanabe's work. It governs the rate at which $X$ makes excursions from a cell corresponding to equality of some $n = k + l$ coordinates into a cell that corresponds to two sets of equalities: one between $k$ coordinates and the other between $l$ coordinates. If we think of $X$ as specifying the location of $N$ particles in $\mathbf{R}$, then in a certain excursion theoretic sense, $\theta(k\!:\!l)$ describes the rate at which a group $n$ particles moving together splits into two groups one consisting of $k$ particles and the other of $l$ particles.

A stochastic flow of kernels on a measurable space $(E, \mathcal{E})$ is a doubly indexed family $(K_{s,t}; s \leq t)$ of random $E \times \mathcal{E}$ transition kernels satisfying the flow property

$$(1.3) \qquad K_{s,u}(x, A) = \int_E K_{s,t}(x, dy) K_{t,u}(y, A), \qquad x \in E, \ A \in \mathcal{E}$$

almost surely for each $s \leq t$. We also postulate independent and stationary innovations in that $K_{t_1,t_2}, K_{t_2,t_3}, \ldots, K_{t_{n-1},t_n}$ are independent for all choices of $t_1 < t_2 < \cdots < t_n$ and $K_{s,t} \stackrel{\text{dist}}{=} K_{s+h,t+h}$ for all $s < t$ and $h$. The general theory of such flows was developed by Le Jan and Raimond in [10]; see also Tsirelson [16]. Stochastic flows of kernels have arisen in the study of the Kraichnan model for turbulent flows; see Falkovich, Gawędzki and Vergassola [3]. Indeed, flows whose two point motions are sticky have been obtained in this context; see Gawędzki and Horvai [4].

One possible interpretation for the flow $K$ is as describing the random evolution of a distribution of mass on $E$. In this case $K_{s,t}(x, A)$ represents the proportion of that mass which was located at $x$ at time $s$ which is within the set $A$ at time $t$. An alternative interpretation of $K$ is as a random environment—in time and space—governing the motion of a particle. Then $K_{s,t}(x, A)$ is the conditional probability given the environment that a particle which is located at $x$ at time $s$ is located within the set $A$ at time $t$.

Taking the space $E$ to be $\mathbf{Z}$, the integer lattice, we may construct a simple example of a stochastic flow of kernels as follows. Fix a probability distribution $\mu$ on $[0, 1]$. Let $\Lambda$ be a Poisson point process on $\mathbf{R} \times \mathbf{Z}$ with uniform intensity, and attach to each point of $(t, x) \in \Lambda$ an independent random variable $R(t, x)$ having the distribution $\mu$. Now consider a particle moving on $\mathbf{Z}$ whose motion is governed by the environment $(R(t, x); (t, x) \in \Lambda)$ as follows. The trajectory of the particle jumps at, and only at, space–time points $(t, x) \in \Lambda$, and the jump is from $x$ to $x + 1$ with probability $R(t, x)$, and from $x$ to $x - 1$ with probability $1 - R(t, x)$. Defining $K_{s,t}(x, A)$ to be the conditional probability given the environment that the particle, when started located at $x$ at time $s$, is located within the set $A$ at time $t$, we obtain a flow of kernels.

Suppose that we have a sequence of probability measures $(\mu_n; n \geq 1)$, centered in that $\int_0^1 x \mu_n(dx) = 1/2$, and with the property that, as $n$ tends to infinity,

$$(1.4) \qquad \sqrt{n} x(1-x) \mu_n(dx) \qquad \text{converges weakly to } \nu(dx),$$

where $\nu$ is some finite measure on $[0,1]$. By weak convergence of finite measures, we mean weak-$*$ convergence induced by the dual space of bounded continuous functions on $[0,1]$. Let $K^n$ be the flow of kernels on $\mathbf{Z}$, associated with $\mu_n$ in the manner described above. It is proved in Section 8, Theorem 8.1, that as $n$ tends to infinity the sequence of flows $K^n$, suitably scaled, converges in law to a flow $K$ on the real line whose $N$-point motions form a consistent family of Brownian motions. Moreover the parameters $\theta(k\!:\!l)$ which specify the corresponding martingale problems satisfy

$$(1.5) \qquad \theta(k\!:\!l) = \int_0^1 x^{k-1}(1-x)^{l-1} \nu(dx) \qquad \text{for } k,l \geq 1.$$

The main difficulties arising in this paper relate to the singular nature of the generator $\mathcal{A}_N^\theta$ for our $N$-dimensional diffusion associated with the family of parameters $\theta$. The action of $\mathcal{A}_N^\theta$ on $C^2$ functions does not characterize the process; in fact, this action does not involve any of the $\theta$ parameters. This means we are unable to appeal to the theory described by Stroock and Varadhan [14]. Instead, our martingale problem is based on the action of generator on a certain vector space of piecewise linear functions. Uniqueness for solutions of the martingale problem reduces, by induction on the dimension $N$ and localization, to showing that the way in which the process leaves the diagonal of $\mathbf{R}^N$ has been specified uniquely. In particular, we study the exit distribution of the process from a $\varepsilon$-neighborhood of the diagonal and determine how, asymptotically as $\varepsilon$ tends to zero, it is described in terms of the $\theta$ parameters. Existence is shown by considering the scaling limit of some Markov chains on the integer lattice $\mathbf{Z}^N$. In fact, these are the $N$-point motions associated with the flows of kernels on $\mathbf{Z}$ described above. Once again the singularity of the generator is an issue. We will be considering functions $f$ for which $\mathcal{A}_N^\theta f$ is not continuous, only upper semi-continuous. Consequently weak convergence will only give us supermartingales. Luckily, this turns out to be enough.

This paper is arranged as follows. In Section 2 we describe the martingale problems that we use to specify our multidimensional diffusions. Sections 3 through 5 are given over to proving uniqueness for these martingale problems, and Section 6 contains a proof of existence. Section 7 contains some auxiliary results that are needed in the proof of uniqueness. Section 8 is concerned with the associated flows of kernels.

Finally let us mention some recent, connected work. In [6] we show that the flow of kernels corresponding to taking the measure $\nu$ equal to a multiple of $\delta_0 + \delta_1$ can be constructed by a filtering procedure applied to a coupled pair of Brownian webs. In [15], Sun and Swart construct a new object, the Brownian net, which we believe is closely related to the flow of kernels in which $\nu$ is a multiple of $\delta_{1/2}$.

**2. A martingale problem.** Consider a pair of $\theta$-coupled Brownian motions $(X_1, X_2)$ as defined in the Introduction. Using Tanaka's formula we can re-express (1.2) as

$$(2.1) \quad |X_1(t) - X_2(t)| - 2\theta \int_0^t \mathbf{1}_{(X_1(s) = X_2(s))} \, ds \qquad \text{is a martingale.}$$

We wish to develop a multidimensional generalization of this statement.

We begin by partitioning $\mathbf{R}^N$ into cells. A cell $E \subset \mathbf{R}^N$ is determined by some weak total ordering $\preceq$ of the $\{1, 2, \ldots, N\}$ via

$$(2.2) \qquad E = \{x \in \mathbf{R}^N : x_i \leq x_j \text{ if and only if } i \preceq j\}.$$

Thus $\{x \in \mathbf{R}^3 : x_1 = x_2 = x_3\}$, $\{x \in \mathbf{R}^3 : x_1 < x_2 = x_3\}$ and $\{x \in \mathbf{R}^3 : x_1 > x_2 > x_3\}$ are three of the thirteen distinct cells into which $\mathbf{R}^3$ is partitioned.

Suppose that $I$ and $J$ are disjoint subsets of $\{1, 2, \ldots, N\}$ with not both $I$ and $J$ empty. With such a pair we associate a vector $v = v_{IJ}$ belonging to $\mathbf{R}^N$ with components given by

$$(2.3) \qquad v_i = \begin{cases} 0, & \text{if } i \notin I \cup J, \\ +1, & \text{if } i \in I, \\ -1, & \text{if } i \in J. \end{cases}$$

We want to associate with each point $x \in \mathbf{R}^N$ certain vectors of this form. To this end, note that each point $x \in \mathbf{R}^N$ determines a partition $\pi(x)$ of $\{1, 2, \ldots, N\}$ such that $i$ and $j$ belong to the same component of $\pi(x)$ if and only if $x_i = x_j$. Then to each point $x \in \mathbf{R}^N$ we associate the set of vectors, denoted by $\mathcal{V}(x)$, which consists of every vector of the form $v = v_{IJ}$ where $I \cup J$ forms one component of the partition $\pi(x)$. The geometric role of the vectors $\mathcal{V}(x)$ should be understood in relation to cells. If $v_{IJ} \in \mathcal{V}(x)$ is such that either $I$ or $J$ is empty then $v_{IJ}$ points from $x$ in a direction that remains in the cell $E$ to which $x$ belongs. We will write $\mathcal{V}_0(x)$ for the subset of $\mathcal{V}(x)$ containing such $v$. If on the other hand both $I$ and $J$ are nonempty then $v_{IJ}$ points from $x$ into another cell which we then call a neighbor of the cell $E$ which contains $x$. We will write $\mathcal{V}_+(x)$ for the subset of $\mathcal{V}(x)$ containing such $v$.

Let $L_N$ be the space of real-valued functions defined on $\mathbf{R}^N$ which are continuous, and whose restriction to each cell is given by a linear function. Given a set of parameters $(\theta(k\!:\!l); k, l \geq 0)$ we define the operator $\mathcal{A}_N^\theta$ from $L_N$ to the space of real valued functions on $\mathbf{R}^N$ which are constant on each cell by

$$(2.4) \qquad \mathcal{A}_N^\theta f(x) = \sum_{v \in \mathcal{V}(x)} \theta(v) \nabla_v f(x).$$

Here on the right-hand side $\theta(v) = \theta(k:l)$ where $k = |I|$ is the number of elements in $I$ and $l = |J|$ is the number of elements in $J$ for $I$ and $J$ determined by $v = v_{IJ}$. $\nabla_v f(x)$ denotes the (one-sided) gradient of $f$ in the direction $v$ at the point $x$, that is,

$$(2.5) \qquad \nabla_v f(x) = \lim_{\varepsilon \downarrow 0} \frac{1}{\varepsilon}(f(x + \varepsilon v) - f(x)).$$

Notice that $\theta(0:0)$ plays no part in the definition of $\mathcal{A}_N^\theta$. Let us also remark that if $\theta$ and $\tilde\theta$ are two parameter families, satisfying for some $\alpha \in \mathbf{R}$,

$$(2.6) \qquad \tilde\theta(k:l) = \theta(k:l) + \alpha \mathbf{1}_{(k=0)} + \alpha \mathbf{1}_{(l=0)},$$

then by the linearity of $f$ within cells, $\mathcal{A}_N^\theta = \mathcal{A}_N^{\tilde\theta}$.

DEFINITION 2.1. We say a continuous, $\mathbf{R}^N$-valued stochastic process $(X(t); t \geq 0)$ solves the $\mathcal{A}_N^\theta$-martingale problem if for each $f \in L_N$,

$$f(X(t)) - \int_0^t \mathcal{A}_N^\theta f(X(s)) \, ds \qquad \text{is a martingale,}$$

relative to some common filtration, and the bracket between coordinates $X_i$ and $X_j$ is given by

$$\langle X_i, X_j \rangle(t) = \int_0^t \mathbf{1}(X_i(s) = X_j(s)) \, ds \qquad \text{for } t \geq 0.$$

In particular $\langle X_i \rangle(t) = t$.

We will sometimes also refer to a probability measure $\mathbf{P}$ on the space of paths $\mathbf{C}([0, \infty), \mathbf{R}^N)$ as being a solution to the $\mathcal{A}_N^\theta$-martingale problem. This means that $X$, being the coordinate process on this space, governed by $\mathbf{P}$, is a solution to $\mathcal{A}_N^\theta$-martingale problem in the sense of the preceeding definition.

The vector space $L_N$ is finite dimensional and the $\mathcal{A}_N^\theta$-martingale problem requires only a finite number of processes to be martingales, but additionally specifies the brackets $\langle X_i, X_j \rangle$. In the case $N = 1$, $L_N$ is one-dimensional and contains only constant multiples of the function $f(x) = x$. The $\mathcal{A}_N^\theta$-martingale problem reduces to the requirement that $X(t) - \beta t$ is a martingale, for $\beta = \theta(1:0) - \theta(0:1)$, and that $\langle X \rangle(t) = t$. Thus, by Lévy's characterization of Brownian motion, the process $(X(t); t \geq 0)$ is a Brownian motion with drift $\beta$. If the consistency property (2.7) holds, then by Proposition 2.1 below, for any dimenison $N$, each component $(X_i(t); t \geq 0)$ of a solution to the $\mathcal{A}_N^\theta$-martingale problem will likewise be a Brownian motion with drift $\beta$. In the case $N = 2$, the space $L_N$ is spanned by $\{x_1, x_2, |x_1 - x_2|\}$. For $f(x) = |x_1 - x_2|$ we find that $\mathcal{A}_N^\theta f(x) = 4\theta(1:1)\mathbf{1}(x_1 = x_2)$, and consequently, assuming (2.7) holds, each pair of components $(X_i, X_j)$ of a solution to the $\mathcal{A}_N^\theta$-martingale problem are $\theta$-coupled Brownian motions (with drift $\beta$) for $\theta = 2\theta(1:1)$. For higher values of $N$, the space $L_N$ is not spanned by the functions $x_i$ and $|x_i - x_j|$; additional functions are used to specify the interactions between more than two particles.

PROPOSITION 2.1. *Suppose that $\theta$ possesses the consistency property*

(2.7) $\qquad \theta(k:l) = \theta(k+1:l) + \theta(k:l+1) \qquad \text{for all } k, l \geq 0.$

*Suppose that $X$ is a solution to the $\mathcal{A}_N^\theta$-martingale problem, and let $Y$ be the process consisting of some $M$ given coordinates of $X$. Then $Y$ is a solution to the $\mathcal{A}_M^\theta$-martingale problem.*

PROOF. We begin by observing the simple fact that if $\sigma$ is a permutation of $\{1, 2, \ldots, N\}$, then $(X_\sigma(t); t \geq 0)$ solves the $\mathcal{A}_N^\theta$-martingale problem if $(X(t); t \geq 0)$ does. Consequently it is enough to prove the proposition in the case that $Y$ consists of the first $N - 1$ coordinates of $X$.

Define $\rho : \mathbf{R}^N \to \mathbf{R}^{N-1}$ to be the projection onto the first $N-1$ coordinates. Suppose that $g \in L_{N-1}$, and let $f = g \circ \rho$, which belongs to $L_N$. $X$ being a solution to the $\mathcal{A}_N^\theta$-martingale problem implies that

$$f(X(t)) - \int_0^t \mathcal{A}_N^\theta f(X(s))\, ds \qquad \text{is a martingale}$$

and since $f(X(t)) = g(Y(t))$ we need to show that $\mathcal{A}_N^\theta f(X(s)) = \mathcal{A}_{N-1}^\theta g(Y(s))$. For this we verify that $\mathcal{A}_N^\theta f(x) = \mathcal{A}_{N-1}^\theta g(\rho(x))$ for all $x \in \mathbf{R}^N$. Fix some $x \in \mathbf{R}^N$, and let $y = \rho(x)$. We will show that

(2.8) $\qquad \displaystyle\sum_{v \in \mathcal{V}(x)} \theta(v) \nabla_v f(x) = \sum_{u \in \mathcal{V}(y)} \theta(u) \nabla_u g(y),$

by matching terms on the right-hand side to terms on the left-hand side.

With this in mind, let us first observe that whenever $I$ and $J$ are disjoint subsets of $\{1, 2, \ldots, N\}$ with $\bar{I} = I \cap \{1, 2, \ldots, N-1\}$ and $\bar{J} = J \cap \{1, 2, \ldots, N-1\}$ not both empty, we have $\rho(v_{IJ}) = u_{\bar{I}\bar{J}}$, where $u = u_{\bar{I}\bar{J}}$ is a vector in $\mathbf{R}^{N-1}$ defined in an analogous manner to (2.3). This together with $f = g \circ \rho$ implies that

$$\nabla_{v_{IJ}} f(x) = \nabla_{u_{\bar{I}\bar{J}}} g(y).$$

Recall that we associate with $x$ a partition $\pi(x)$ of $\{1, 2, \ldots, N\}$ into classes $C_1, C_2, \ldots, C_n$ say, without loss of generality, the class $C_1$ containing $N$. Then $\pi(y)$ is the partition of $\{1, 2, \ldots, N-1\}$ into classes $C_1 \setminus \{N\}, C_2, \ldots, C_n$. It may be that $C_1$ is just the singleton $\{N\}$ in which case we have only classes $C_2, \ldots, C_n$ here. Each term on the left-hand side of (2.8) corresponds to splitting some class $C_r$ of $\pi(x)$ into an ordered pair $(I, J)$ of parts, one of which may be empty. There are two different cases to consider for this pair.

First, when $I \cup J = C_r$ for $r \neq 1$ which consequently does not contain $N$. Then by taking $\bar{I} = I$ and $\bar{J} = J$ we obtain a corresponding term on the right-hand side. Moreover we have

$$\theta(v_{IJ})\nabla_{v_{IJ}}f(x) = \theta(u_{\bar{I}\bar{J}})\nabla_{u_{\bar{I}\bar{J}}}g(y)$$

and the contributions to the two sides of (2.8) are equal.

In the second case, we consider $I \cup J = C_1$. Without loss of generality suppose that $N \in I$. Then we consider a second pair $(I', J')$ of subsets of $\{1, 2, \ldots, N\}$ with $I' = I \setminus \{N\}$ and $J' = J \cup \{N\}$. Set $\bar{I} = I'$ and $\bar{J} = J$. If we suppose that there are $k$ elements in $\bar{I}$ and $l$ elements in $\bar{J}$, then, by the consistency property for $\theta$,

$$\theta(v_{IJ}) + \theta(v_{I'J'}) = \theta(k+1:l) + \theta(k:l+1) = \theta(k:l) = \theta(u_{\bar{I}\bar{J}}).$$

Thus, provided not both $\bar{I}$ and $\bar{J}$ are empty the net contribution of $v_{IJ}$ and $v_{I'J'}$ to the left-hand side of (2.8) matches the contribution of $u_{\bar{I}\bar{J}}$ to the right-hand side. Actually, even in the case $\bar{I}$ and $\bar{J}$ are both empty, this is still true. Then there is no term on the right-hand side, but $f$ does not depend on $x_N$ so we have $\nabla_{v_{IJ}}f(x) = \nabla_{v_{I'J'}}f(x) = 0$.

Taken together these two cases for $(I, J)$ exhaust all possibilities, and we have also accounted for all terms on the right-hand side of (2.8). The equality $\mathcal{A}_N^\theta f(x) = \mathcal{A}_{N-1}^\theta g(y)$ is thus proven.

The proof of the proposition is completeted by noting that the brackets $\langle Y_i, Y_j \rangle$ are as required. $\square$

The role of the parameters $\theta(k:l)$ with both $k$ and $l$ strictly positive was explained in the Introduction. They may be interpreted loosely as rates and as such we impose the positivity condition

(2.9) $\qquad \theta(k:l) \geq 0 \qquad \text{for all } k, l \geq 1.$

The parameters $\theta(k:0)$ and $\theta(0:l)$ are not necessarily positive. Their role is probably best described as contributing correction terms to the generator $\mathcal{A}_N^\theta$ which ensure the consistency as $N$ varies. As noted above, if $X$ solves the $\mathcal{A}_N^\theta$-martingale problem for a family of parameters $\theta$ satisfying (2.7), then each coordinate $(X_i(t); t \geq 0)$ is a Brownian motion with drift $\beta = \theta(1:0) - \theta(0:1)$. Observe that any consistent family $(\theta(k:l); k, l \geq 0)$ is determined completely by the restricted family $(\theta(k:l); k, l \geq 1)$ of "splitting rates" together with the value of the drift $\beta = \theta(1:0) - \theta(0:1)$. The following proposition shows that we can eliminate drift in the usual way.

PROPOSITION 2.2. *Suppose that $X$ solves the $\mathcal{A}_N^\theta$-martingale problem and let $\tilde{X}_i(t) = X_i(t) - 2\beta t$. Then $\tilde{X}$ solves the $\mathcal{A}_N^{\tilde{\theta}}$-martingale problem where*

(2.10) $$\tilde{\theta}(k:l) = \theta(k:l) - \beta \mathbf{1}_{(k=0)} + \beta \mathbf{1}_{(l=0)}.$$

PROOF. Let $\mathbf{1}$ denote the vector $(1, 1, \ldots, 1) \in \mathbf{R}^N$. We have, for any $f \in L_N$,

$$f(\tilde{X}(t)) - \int_0^t \mathcal{A}_N^{\tilde{\theta}} f(\tilde{X}(s))\,ds = f(X(t)) - 2\beta t f(\mathbf{1}) - \int_0^t \mathcal{A}_N^{\tilde{\theta}} f(X(s))\,ds,$$

using the facts that $f$ is linear when restricted to each cell, and $\mathcal{A}_N^{\tilde{\theta}} f$ is constant on each cell. Thus, it is enough to verify that $2\beta f(\mathbf{1}) + \mathcal{A}_N^{\tilde{\theta}} f(x) = \mathcal{A}_N^\theta f(x)$ for all $x \in \mathbf{R}^N$. Since $\tilde{\theta}(k:l) = \theta(k:l)$ whenever $k, l \geq 1$, this comes down to checking that

$$2\beta f(\mathbf{1}) + \sum_{v \in \mathcal{V}_0(x)} \tilde{\theta}(v) \nabla_v f(x) = \sum_{v \in \mathcal{V}_0(x)} \theta(v) \nabla_v f(x).$$

This holds by virtue of the linearity of $f$ within the cell containing $x$. □

The main purpose of this paper is to prove the following theorem.

THEOREM 2.1. *Let $\theta$ be a family of parameters satisfying the consistency and positivity properties given above in (2.7) and (2.9). For each $N \geq 1$ and $x \in \mathbf{R}^N$ there exists a process solving the $\mathcal{A}_N^\theta$-martingale problem starting from $x$. Moreover the law of this process is unique.*

**3. Leaving the diagonal.** Throughout this section we assume that we have fixed some family $\theta$ of parameters satisfying the consistency and positivity conditions given in (2.7) and (2.9). We suppose that $X$ is some solution to the $\mathcal{A}_N^\theta$-martingale problem. Let

$$T_\varepsilon = \inf\{t \geq 0 : |X_i(t) - X_j(t)| \geq \varepsilon \text{ for some } i, j \in \{1, 2, \ldots, N\}\}.$$

We intend show that if $X$ starts on the diagonal $D = \{x \in \mathbf{R}^N : x_i = x_j \text{ for all } i, j\}$ then the distribution of $X(T_\varepsilon)$ is, for small $\varepsilon$, concentrated on the set of cells $E \subset \mathbf{R}^N$ which are neighbors of $D$. Moreover the parameters $\theta(k:l)$ determine the way the exit distribution is apportioned between these cells. To be more precise, let $\mathcal{V}_+(D)$ denote the collection of vectors of the form $v_{IJ}$ where $I \cup J = \{1, 2, \ldots, N\}$ and neither $I$ nor $J$ is empty. Such a vector determines a cell

(3.1) $E(v) = \{y \in \mathbf{R}^N : y = x + \beta v \text{ for } x \in D \text{ and } \beta \in \mathbf{R} \text{ with } \beta > 0\}.$

These cells are the neighbors of $D$.

The main result of this section is as follows.

THEOREM 3.1. *Let $X$ start from the diagonal. Then the following limits exist and are determined by the family of parameters $\theta$:*

$$\lim_{\varepsilon \to 0} \frac{1}{\varepsilon} \mathbf{E}[T_\varepsilon] = \frac{1}{2 \sum_{v \in \mathcal{V}_+(D)} \theta(v)}$$

*and for each cell $E(v)$, a neighbor of $D$,*

$$\lim_{\varepsilon \to 0} \mathbf{P}(X(T_\varepsilon) \in E(v)) = \frac{\theta(v)}{\sum_{u \in \mathcal{V}_+(D)} \theta(u)}.$$

If the dimension $N = 2$, then the result is easily established, and in the sequel we shall mostly assume $N \geq 3$. We notice also that by virtue of Proposition 2.2 and the remark made at (2.6) we can make the simplifying assumption that $\theta(0\!:\!1) = \theta(1\!:\!0) = 0$, and so in particular every coordinate of $X$ is a driftless Brownian motion.

LEMMA 3.1. *Let $X$ start from the diagonal then*

$$\mathbf{E}\left[\int_0^{T_\varepsilon} \mathbf{1}(X(s) \notin D)\, ds\right] \leq \frac{N(N-1)}{4} \varepsilon^2.$$

PROOF. Each pair of coordinates $(X_i(t), X_j(t))$ is, by Proposition 2.1, a pair of $\theta$-coupled Brownian motions with $\theta = 2\theta(1\!:\!1)$. Consequently

$$(X_i(t) - X_j(t))^2 - 2 \int_0^t \mathbf{1}(X_i(s) \neq X_j(s))\, ds$$

is a martingale.

Applying the optional stopping theorem at time $T_\varepsilon \wedge t$, and letting $t \to \infty$, we obtain,

$$\mathbf{E}\left[\int_0^{T_\varepsilon} \mathbf{1}(X_i(s) \neq X_j(s))\, ds\right] = \frac{1}{2} \mathbf{E}[(X_i(T_\varepsilon) - X_j(T_\varepsilon))^2] \leq \varepsilon^2/2.$$

The result follows by summing over all possible pairs $i$ and $j$. □

In the following series of lemmas we will assume that the dimension $N = 3$. Later we will apply our conclusions to each triple of coordinates of a process $X$ having arbitrary dimension. When $N = 3$ the diagonal $\{x_1 = x_2 = x_3\}$ has 6 neighboring cells each of the form $E_{ij}^k = \{x_i = x_j < x_k\}$ or $E_k^{ij} = \{x_k < x_i = x_j\}$.

Define the process $(Y(t); t \geq 0)$ by

(3.2) $$Y(t) = \inf_{1 \leq i < j \leq 3} |X_i(t) - X_j(t)|.$$

$Y(t)$ measures the distance from $X(t)$ to the union of $D$ and the six cells neighboring $D$. We begin by showing $Y(T_\varepsilon)$ is typically small, and will subsequently show that it is in fact zero with high probability. We denote $2\theta(1\!:\!1)$ by $\theta$ for the rest of this section.

LEMMA 3.2. *Suppose that $X$ starts from the diagonal of $\mathbf{R}^3$ then*

$$\mathbf{E}[Y(T_\varepsilon)] \leq 6\theta\varepsilon^2.$$

PROOF. Each pair of coordinates $(X_i(t), X_j(t))$ is a pair of $\theta$-coupled Brownian motions with $\theta = 2\theta(1\!:\!1)$. Consequently

$$|X_i(t) - X_j(t)| = \int_0^t \text{sgn}(X_i(s) - X_j(s))\, d(X_i(s) - X_j(s))$$
$$+ 2\theta \int_0^t \mathbf{1}_{(X_i(s) = X_j(s))}\, ds.$$

Using this, together with the brackets $\langle X_i, X_j \rangle(t) = \int_0^t \mathbf{1}_{(X_i(s) = X_j(s))}\, ds$, we calculate that

$$Z(t) = \prod_{i<j} |X_i(t) - X_j(t)|$$
$$= M(t) + 2\theta \sum_{i=1}^3 \int_0^t |X_i(s) - X_j(s)||X_i(s) - X_k(s)| \mathbf{1}_{(X_j(s) = X_k(s) \neq X_i(s))}\, ds,$$

where in the last expression $\{j,k\} = \{1,2,3\} \setminus \{i\}$, and where $M(t)$ is given by a stochastic integral with respect to $X(t)$ and is hence a local martingale. Applying the optional stopping theorem at $T_\varepsilon \wedge t$, noting that $M(s)$ is bounded on the interval $[0, T_\varepsilon \wedge t]$, and then letting $t$ tend to infinity, we obtain

$$\mathbf{E}[Z(T_\varepsilon)] = 2\theta \mathbf{E}\left[\sum_{i=1}^3 \int_0^{T_\varepsilon} |X_i(s) - X_j(s)||X_i(s) - X_k(s)| \mathbf{1}_{(X_j(s) = X_k(s) \neq X_i(s))}\, ds\right]$$
$$\leq 2\theta\varepsilon^2 \mathbf{E}\left[\int_0^{T_\varepsilon} \mathbf{1}(X(s) \notin D)\, ds\right] \leq 3\theta\varepsilon^4,$$

appealing to Lemma 3.1 for the last inequality. To conclude we note that

$$Z(T_\varepsilon) = \varepsilon Y(T_\varepsilon)(\varepsilon - Y(T_\varepsilon)) \geq \frac{\varepsilon^2}{2} Y(T_\varepsilon). \qquad \square$$

LEMMA 3.3. *Suppose that $X$ starts from a point $x \in \mathbf{R}^3$ satisfying $\sup_{i \leq j} |x_i - x_j| = \varepsilon$. Let $T_{0,2\varepsilon}$ be the first exit time of $X$ from the domain $\{x \in \mathbf{R}^3 : 0 < \sup_{i \leq j} |x_i - x_j| < 2\varepsilon\}$. Then*

$$\mathbf{P}(Y(T_{0,2\varepsilon}) \neq 0) \leq C(y\varepsilon^{-1} + \theta\varepsilon),$$

*where $y = \inf_{i<j} |x_i - x_j|$, and $C$ is some universal constant.*

PROOF. By permuting coordinates we may suppose that $y = |x_1 - x_2| \leq \min(|x_1 - x_3|, |x_2 - x_3|)$. We will also assume that $x_1, x_2 > x_3$, the other possibility that $x_1, x_2 < x_3$ can be dealt with by a analogous argument. Let

$$U^\varepsilon = \{x \in \mathbf{R}^3 : x_1, x_2 > x_3\} \cap \left\{x \in \mathbf{R}^3 : 0 < \sup_{i \leq j} |x_i - x_j| < 2\varepsilon\right\},$$

and $T(U^\varepsilon)$ denote the first time $X$ exits this domain. Let $V(t) = \frac{1}{\sqrt{2}}|X_1(t) - X_2(t)|$ and $W(t) = \frac{1}{\sqrt{6}}(X_1(t) + X_2(t) - 2X_3(t))$. Then the following equality in law holds between stopped processes:

$$(W(t \wedge T(U^\varepsilon)), V(t \wedge T(U^\varepsilon)); t \geq 0) \stackrel{\text{law}}{=} (\xi(t \wedge \tau), \eta(t \wedge \tau); t \geq 0),$$

where $(\eta, \xi)$ is a sticky Brownian motion in the half plane of the type described in Section 7 with data $(a_0, \theta_0) = (4/3, \sqrt{2}\theta)$ and $\tau$ is the first exit time of $(\xi, \eta)$ from the triangle

$$\triangle(\varepsilon) = \{(x, y) \in \mathbf{R}^2 : y \geq 0, \sqrt{6}x - \sqrt{2}y > 0, \sqrt{6}x + \sqrt{2}y < 4\varepsilon\}.$$

The desired inequality holds by noting that the event $\{Y(T_{0,2\varepsilon}) \neq 0\}$ is contained in the event $\{V(T(U^\varepsilon)) \neq 0\}$ and applying Proposition 7.2. □

LEMMA 3.4. *Suppose that $X$ starts from the diagonal of $\mathbf{R}^3$. Then*

$$\mathbf{P}(Y(T_{2\varepsilon}) \neq 0) \leq C\theta\varepsilon,$$

*where $C$ is some universal constant.*

PROOF. Let $R(t) = \sup_{i \leq j} |X_i(t) - X_j(t)|$. Define a sequence of stopping times by $T_0 = 0$, and then in general,

$$T_{2n+1} = \inf\{t \geq T_{2n} : R(t) = \varepsilon\},$$
$$T_{2n} = \inf\{t \geq T_{2n-1} : R(t) = 0\}.$$

It is easy to see that these stopping times are almost surely finite. Denote the natural filtration of $X$ by $(\mathcal{F}_t; t \geq 0)$. Then $R(t)$ is a $\mathcal{F}_t$-submartingale and consequently,

$$2\varepsilon \mathbf{P}(T_{2n} > T_{2\varepsilon} > T_{2n-1}|\mathcal{F}_{T_{2n-1}})$$
$$= \mathbf{E}[R(T_{2\varepsilon} \wedge T_{2n})\mathbf{1}_{(T_{2\varepsilon} > T_{2n-1})}|\mathcal{F}_{T_{2n-1}}]$$
$$\geq R(T_{2n-1})\mathbf{1}_{(T_{2\varepsilon} > T_{2n-1})} = \varepsilon \mathbf{1}_{(T_{2\varepsilon} > T_{2n-1})}.$$

From this it follows that

$$\mathbf{P}(T_{2\varepsilon} > T_{2n+1}) \leq \tfrac{1}{2}\mathbf{P}(T_{2\varepsilon} > T_{2n-1})$$

and so $\mathbf{P}(T_{2\varepsilon} > T_{2n}) = \mathbf{P}(T_{2\varepsilon} > T_{2n+1}) \leq (\tfrac{1}{2})^n$.

By well-known properties of martingale problems (see Stroock and Varadhan [14]), the conditional law of the process $(X(T_{2n-1} + t); t \geq 0)$ given $\mathcal{F}_{T_{2n-1}}$ is almost surely a solution to the $\mathcal{A}_N^\theta$-martingale problem starting from $X(T_{2n-1})$. Thus a conditional version of Lemma 3.3 holds

$$\mathbf{P}(Y(T_{0,2\varepsilon}^n) \neq 0 | \mathcal{F}_{T_{2n+1}}) \leq C(Y(T_{2n+1})\varepsilon^{-1} + \theta\varepsilon) \quad \text{a.s.},$$

where $T_{0,2\varepsilon}^n = \inf\{t \geq T_{2n+1} : R(t) = 0 \text{ or } 2\varepsilon\}$. Similarly the following conditional version of Lemma 3.2 holds:

$$\mathbf{E}[Y(T_{2n+1})|\mathcal{F}_{T_{2n}}] \leq 6\theta\varepsilon^2 \quad \text{a.s.}$$

Combining these two estimates gives

$$\mathbf{P}(Y(T_{0,2\varepsilon}^n) \neq 0 | \mathcal{F}_{T_{2n}}) \leq 7C\theta\varepsilon \quad \text{a.s.}$$

Finally, turning to the quantity of interest

$$\mathbf{P}(Y(T_{2\varepsilon}) \neq 0)$$
$$= \sum_{n=0}^{\infty} \mathbf{P}(Y(T_{2\varepsilon}) \neq 0 \text{ and } T_{2n+1} < T_{2\varepsilon} < T_{2n+2})$$
$$= \sum_{n=0}^{\infty} \mathbf{P}(Y(T_{0,2\varepsilon}^n) \neq 0 \text{ and } T_{2n+1} < T_{2\varepsilon} < T_{2n+2})$$
$$\leq \sum_{n=0}^{\infty} \mathbf{P}(Y(T_{0,2\varepsilon}^n) \neq 0 \text{ and } T_{2n} < T_{2\varepsilon}) \leq \sum_{n=0}^{\infty} 7C\theta\varepsilon \left(\frac{1}{2}\right)^n = 14C\theta\varepsilon. \quad \square$$

We now return to considering $X$ having arbitrary dimension $N \geq 3$.

PROPOSITION 3.1. *Let $X$ start from the diagonal of $\mathbf{R}^N$. Let $\Lambda$ be the event that there are three or more distinct values taken by the coordinates of $X(T_\varepsilon)$. Then*

$$\mathbf{P}(\Lambda) \leq CN^3\theta\varepsilon,$$

*where $C$ is a universal constant.*

PROOF. For each triple of distinct indices $(i, j, k)$ define

$$T_{ijk} = \inf\{t \geq 0 : \max(X_i(t), X_j(t), X_k(t)) - \min(X_i(t), X_j(t), X_k(t)) = \varepsilon\}.$$

The graph of $T_\varepsilon$ is contained in the union of the graphs of $T_{ijk}$ as $(i, j, k)$ varies and

$$\Lambda \subseteq \bigcup_{i,j,k} \Lambda_{ijk},$$

where $\Lambda_{ijk}$ is the event that $X_i(T_{ijk}), X_j(T_{ijk})$ and $X_k(T_{ijk})$ are distinct. By virtue of Proposition 2.1 the process $(X_i(t), X_j(t), X_k(t); t \geq 0)$ solves the $\mathcal{A}_3^\theta$-martingale problem starting from the diagonal, and consequently we deduce from Lemma 3.4, with $2\varepsilon$ replaced by $\varepsilon$ in its statement, that $\mathbf{P}(\Lambda_{ijk}) \leq C\theta\varepsilon$. The result follows by summing over $i, j$ and $k$. □

PROOF OF THEOREM 3.1. Take $f \in L_N$ to be given by

$$f(x) = \max_{1 \leq i \leq N}\{x_i\} - \min_{1 \leq i \leq N}\{x_i\}.$$

Then it is easily verified that $\nabla_v f(x) = 2$ for $x \in D$ and $v \in \mathcal{V}_+(D)$, while $\nabla_v f(x) = 0$ for $x \in D$ and $v \in \mathcal{V}(x) \setminus \mathcal{V}_+(D)$. Also $|\nabla_v f(x)| \leq 2$ for any $x \in \mathbf{R}^N$ and $v \in \mathcal{V}(x)$. This latter fact implies that $|\mathcal{A}_N^\theta f(x)|$ is bounded by some constant depending on $N$ and $\theta$ for all $x \in \mathbf{R}^N$.

Applying the optional stopping theorem to the martingale $f(X(t)) - \int_0^t \mathcal{A}_N^\theta \times f(X(s))\,ds$ at the time $T_\varepsilon \wedge t$ and letting $t$ tend to infinity we obtain

$$\begin{aligned}
\varepsilon &= \mathbf{E}[f(X(T_\varepsilon))] \\
&= \mathbf{E}\left[\int_0^{T_\varepsilon} \mathcal{A}_N^\theta f(X(s))\,ds\right] \\
&= \mathbf{E}\left[\int_0^{T_\varepsilon} \mathcal{A}_N^\theta f(X(s))\mathbf{1}(X(s) \notin D)\,ds\right] \\
&\quad + \sum_{v \in \mathcal{V}_+(D)} 2\theta(v)\mathbf{E}\left[\int_0^{T_\varepsilon} \mathbf{1}(X(s) \in D)\,ds\right].
\end{aligned}$$

The first assertion of the theorem now follows by dividing through by $\varepsilon$ and passing to the limit with an appeal to Lemma 3.1.

For the second result consider some $E(v)$, a cell which neighbors $D$. Then for some partition of $\{1, 2, \ldots, N\}$ into nonempty disjoint classes $I$ and $J$ we have $v = v_{IJ}$ and

$$E(v) = \{x \in \mathbf{R}^N : x_i > x_j \text{ if } i \in I, j \in J, \text{ and } x_i = x_j \text{ if either } i, j \in I \text{ or } i, j \in J\}.$$

We associate with $E(v)$ a function $f_v \in L_N$, defined with the aid of the partition as

$$f_v(x) = \min_{i \in I, j \in J}(x_i - x_j)^+.$$

Once again it is easy to verify that for any $x \in D$, we have $\nabla_v f_v(x) = 2$, while $\nabla_u f_v(x) = 0$ for any $u \in \mathcal{V}(x)$ with $u \neq v$. Also $|\nabla_u f_v(x)| \leq 2$ for any $x \in \mathbf{R}^N$ and $u \in \mathcal{V}(x)$ and, as a consequence, $|\mathcal{A}_N^\theta f_v(x)|$ is bounded by some constant depending on $N$ and $\theta$ for all $x \in \mathbf{R}^N$. Applying the optional stopping theorem to the martingale $f_v(X(t)) - \int_0^t \mathcal{A}_N^\theta f_v(X(s))\,ds$ at the time $T_\varepsilon \wedge t$ and letting $t$ tend to infinity we obtain

$$\mathbf{E}[f_v(X(T_\varepsilon))] = \mathbf{E}\left[\int_0^{T_\varepsilon} \mathcal{A}_N^\theta f_v(X(s))\,ds\right]$$
$$= \mathbf{E}\left[\int_0^{T_\varepsilon} \mathcal{A}_N^\theta f_v(X(s))\mathbf{1}(X(s) \notin D)\,ds\right]$$
$$+ 2\theta(v)\mathbf{E}\left[\int_0^{T_\varepsilon} \mathbf{1}(X(s) \in D)\,ds\right].$$

We divide through by $\varepsilon$ and pass to the limit. On the right-hand side, by Lemma 3.1 and the first part of the current theorem, we obtain

$$\frac{\theta(v)}{\sum_{u \in \mathcal{V}_+(D)} \theta(u)}.$$

While on the left-hand side, by virtue of Proposition 3.1, we have

$$\lim_{\varepsilon \to 0} \frac{1}{\varepsilon}\mathbf{E}[f_v(X(T_\varepsilon))] = \lim_{\varepsilon \to 0} \frac{1}{\varepsilon}\mathbf{E}[f_v(X(T_\varepsilon))\mathbf{1}_{\Lambda^c}] = \lim_{\varepsilon \to 0} \mathbf{P}(X(T_\varepsilon) \in E(v)). \quad \square$$

Later it will be important to us to have noted that the proof just given shows not just that the limits appearing in the statement of the theorem hold, but in fact gives us rates of convergence. Indeed we have proved that, for some constant $C$ that depends only on $N$ and $\theta$ we have

$$(3.3) \qquad \left|\mathbf{E}[T_\varepsilon] - \frac{\varepsilon}{2\sum_{v \in \mathcal{V}_+(D)} \theta(v)}\right| \leq C\varepsilon^2 \quad \text{and}$$

$$(3.4) \qquad \left|\mathbf{P}(X(T_\varepsilon) \in E(v)) - \frac{\theta(v)}{\sum_{u \in \mathcal{V}_+(D)} \theta(u)}\right| \leq C\varepsilon.$$

**4. The process stopped at the time of hitting the diagonal.** Consider a partition $\pi = (\pi_1, \pi_2)$ of $\{1, 2, \ldots, N\}$ into two parts. Let $|\pi_i|$ denote the number of elements in the class $\pi_i$. We will say that a $\mathbf{R}^N$-valued process $(X(t); t \geq 0)$ solves the $\mathcal{A}_\pi^\theta$-martingale problem starting from $x$ if, for each $i = 1, 2$, the process $(X_k(t); k \in \pi_i, t \geq 0)$ solves the $\mathcal{A}_{|\pi_i|}^\theta$-martingale problem (relative to the natural filtration of $X$), and if every pair of Brownian motions $(X_k(t); t \geq 0)$ and $(X_l(t); t \geq 0)$ for which $k$ and $l$ belong to distinct parts of $\pi$ are orthogonal.

PROPOSITION 4.1. *Suppose that for each class $\pi_i$ of the partition $\pi$, the $\mathcal{A}_{|\pi_i|}^\theta$-martingale problem has a solution unique in law for and each possible starting point $x \in \mathbf{R}^{|\pi_i|}$. Then the $\mathcal{A}_\pi^\theta$-martingale problem has a solution unique in law for every starting point $x \in \mathbf{R}^N$.*

PROOF. For simplicity we assume $\theta(0\!:\!1) = \theta(1\!:\!0) = 0$, and consider the case $\pi_1 = \{1, 2, \ldots, n\}$ and $\pi_2 = \{n+1, n+2, \ldots, N\}$. Suppose $X$ is a solution to the $\mathcal{A}_\pi^\theta$-martingale problem starting from some point $x \in \mathbf{R}^N$. Let $Y$ be the $\mathbf{R}^n$-valued process consisting of the first $n$ coordinates of $X$, and $Z$ the $\mathbf{R}^{N-n}$-valued process consisting of the last $N-n$ coordinates of $X$. Then $Y$ and $Z$ solve the $\mathcal{A}_{|\pi_1|}^\theta$ and $\mathcal{A}_{|\pi_2|}^\theta$-martingale problems, respectively. The uniqueness-in-law property for these martingale problems implies martingale representation, which we may obtain by applying Theorem 40 of Chapter IV of [12] to the collection of distinct martingales of the form

$$Y_i^E(t) = \int_0^t \mathbf{1}(Y(s) \in E)\, dY_i(s),$$

where $E$ is a cell in $\mathbf{R}^n$ and $i = 1, 2, \ldots, n$ (and similarly for $Z$). As a consequence, for arbitrary $F \in L^2(Y)$ and $G \in L^2(Z)$, we may write

$$F = \mathbf{E}[F] + \int_0^\infty f_s\, dY_s \quad \text{and} \quad G = \mathbf{E}[G] + \int_0^\infty g_s\, dZ_s,$$

for suitable vector-valued predictable processes $(f_s; s \geq 0)$ and $(g_s; s \geq 0)$. But since $Y$ and $Z$ are orthogonal martingales we obtain from this

$$\mathbf{E}[FG] = \mathbf{E}[F]\mathbf{E}[G]$$

and hence deduce that $Y$ and $Z$ are independent. This identifies the law of $X = (Y, Z)$ as being the product of uniquely determined laws, and we are done. $\square$

For each nonempty subset $S$ of $\{1, 2, \ldots, N\}$ having nonempty complement $S^c$, consider the open set

$$U^S = \{x \in \mathbf{R}^N : x_i > x_j \text{ for all } i \in S, j \in S^c\}.$$

For $x \in \mathbf{R}^N$ write $p_S(x)$ for the vector $(x_i; i \in S) \in \mathbf{R}^{|S|}$, and let $p_{S^c}(x)$ be defined similarly. Fix some $S$ and a point $x \in U^S$. Suppose that $X$ solves the $\mathcal{A}_N^\theta$-martingale problem starting from $x$, and let $T_S$ be the first exit time of $X$ from $U^S$. Put $Y(t) = p_S(X(t))$ and $Z(t) = p_{S^c}(X(t))$. By Proposition 2.1 the processes $Y$ and $Z$ solve the $\mathcal{A}_{|S|}^\theta$- and $\mathcal{A}_{|S^c|}^\theta$-martingale problems starting from $y = p_S(x)$ and $z = p_{S^c}(x)$, respectively. Moreover, any pair of coordinates, $(Y_i(T_S \wedge t); t \geq 0)$ and $(Z_j(T_S \wedge t); t \geq 0)$, of the stopped processes are orthogonal. Let $\pi$ denote the partition with parts $S$ and $S^c$. If $T$ is any stopping time almost surely less than or equal to $T_S$, then by standard arguments for stopped martingale problems,

(4.1)  if the solution to the $\mathcal{A}_\pi^\theta$-martingale problem starting from $x$ is known to be unique in law, then the law of $(X(T \wedge t); t \geq 0)$, is also uniquely determined.

Suppose that for every $n \leq N-1$ and $x \in \mathbf{R}^n$ the $\mathcal{A}_n^\theta$-martingale problem starting from $x$ has a solution whose law is uniquely determined. Then by virtue of Proposition 4.1, the hypothesis at (4.1) is verified, and the conclusion that the law of the stopped process is uniquely determined holds. Fix some $\varepsilon > 0$ and observe as $S$ varies the collection of all $U^S$ forms a cover of $K^\varepsilon = \{x \in \mathbf{R}^N : \varepsilon \leq |x_i - x_j| \text{ for some } i \neq j\} \cap \{x \in \mathbf{R}^N : |x_i| \leq 1/\varepsilon \text{ for all } i\}$, having the property that if $x \in K^\varepsilon$ then there exists some subset $S(x)$ so that $U^{S(x)}$ contains the ball $\{y \in \mathbf{R}^N : \|y - x\| < \varepsilon/(2N)\}$. Now fix $x \in K^\varepsilon$, and suppose that $X$ solves the $\mathcal{A}_N^\theta$-martingale problem starting from $x$. Put $T_0 = 0$, and for $i \geq 0$,

(4.2) $$T_{i+1} = \inf\{t \geq T_i : X(t) \notin U^{S(i)} \cap K^\varepsilon\},$$

where $S(i) = S(X(T_i))$. Let $T_\varepsilon$ denote the first time the process $X$ exits $K^\varepsilon$; this is almost surely finite, and consequently by continuity of the paths of $X$, $T_i = T_\varepsilon$, for sufficiently large $i$ with probability one. Now consider the conditional distribution of $(X(T_i + t) \wedge T_{i+1}); t \geq 0)$ given $\mathcal{F}_{T_i}$ where $(\mathcal{F}_t; t \geq 0)$ is the natural filtration of $X$. This conditional law is almost surely a solution to the $\mathcal{A}_N^\theta$-martingale problem starting from $X(T_i)$ and stopped on first exiting $U^{S(i)} \cap K^\varepsilon$. As an application of (4.1) we deduce that this conditional law is uniquely determined. Then it follows, by a standard splicing argument of the type used in [14], that the law of $(X(t \wedge T_\varepsilon); t \geq 0)$ is also unique. Finally letting $\varepsilon$ tend down to zero, we have established the following, in which $T_D$ denotes the first time that $X$ reaches the diagonal $D$ of $\mathbf{R}^N$.

PROPOSITION 4.2. *Suppose that for every $n \leq N-1$ and $x \in \mathbf{R}^n$ the $\mathcal{A}_n^\theta$-martingale problem starting from $x$ has a solution whose law is uniquely determined. Then if $X$ is any solution to the $\mathcal{A}_N^\theta$-martingale problem starting from a point $x \in \mathbf{R}^N$, the law of the stopped process $(X(T_D \wedge t); t \geq 0)$ is uniquely determined.*

**5. The proof of uniqueness.** Fix a family of parameters $\theta$. Anticipating an argument by induction on $N$, we will make the following hypothesis:

(5.1) If $X$ is any solution to the $\mathcal{A}_N^\theta$-martingale problem starting from a point $x \in \mathbf{R}^N$, the law of the stopped process $(X(T_D \wedge t); t \geq 0)$ is uniquely determined.

Now suppose that $X$ is governed by a family of probability measures $(\mathbf{P}_x; x \in \mathbf{R}^N)$, and that under each $\mathbf{P}_x$ it solves the $\mathcal{A}_N^\theta$-martingale problem and starts from $x$. Define the function $\psi_\lambda$ by

(5.2) $$\psi_\lambda(x) = \mathbf{E}_x[\exp(-\lambda T_D)]$$

and similarly for any test function $f$ we define $R_\lambda^0 f$ by

(5.3) $$R_\lambda^0 f(x) = \mathbf{E}_x\left[\int_0^{T_D} e^{-\lambda s} f(X(s))\, ds\right].$$

Our aim is show that the expectation

$$\mathbf{E}_x\left[\int_0^\infty e^{-\lambda s} f(X(s))\,ds\right], \tag{5.4}$$

which a priori depends on possible choices of the family of measures $(\mathbf{P}_x; x \in \mathbf{R}^N)$, is uniquely determined. The hypothesis (5.1) ensures that $R_\lambda^0 f$ and $\psi_\lambda$ do not depend on any such choice. We will show that, whenever $f$ is invariant under translations along the diagonal, and also zero in a neighborhood of it, the expectation (5.4) is given by

$$R_\lambda^0 f(x) + \kappa \psi_\lambda(x), \tag{5.5}$$

for a certain constant $\kappa$ whose value is also unaffected by any freedom in choosing $(\mathbf{P}_x; x \in \mathbf{R}^N)$. Notice that, assuming the uniqueness to the martingale problem that we are intending to prove, and the consequent Markov property, the equality of (5.4) and (5.5) is the usual decomposition of the resolvent $R_\lambda f$ and $\kappa = R_\lambda f(0)$. Moreover, formally applying the generator to this equality gives $\lambda R_\lambda f(0) = \mathcal{A}_N^\theta R_\lambda^0 f(0) + R_\lambda f(0) \mathcal{A}_N^\theta \psi_\lambda(0)$ which gives $R_\lambda f(0)$ in terms of $R_\lambda^0 f$ and $\psi_\lambda$. This motivates the calculations of Propositions 5.1 and 5.2 below, in which $\kappa$ is determined with the help of the description provided by Theorem 3.1 of the way $X$ leaves the diagonal. We need the following lemmas.

LEMMA 5.1. *There exists a constant $C$, depending on $\lambda$, $\theta$ and $N$ alone, such that*

$$1 - \psi_\lambda(x) \leq C\sqrt{\mathrm{dist}(x,D)} \qquad \textit{whenever } \mathrm{dist}(x,D) < 1$$

*and also, for any bounded $f$,*

$$R_\lambda^0 f(x) \leq C\sqrt{\mathrm{dist}(x,D)}\|f\|_\infty \qquad \textit{whenever } \mathrm{dist}(x,D) < 1.$$

PROOF. For $i \neq j$ the difference $(X_i(t) - X_j(t); t \geq 0)$ is distributed as a sticky Brownian motion on $\mathbf{R}$. Using standard results on this process, see, for example, [2], it is easily deduced that

$$\mathbf{P}_x(X_i(t) \neq X_j(t)) \leq \frac{1}{\sqrt{t}}|x_i - x_j| + 4\theta(1:1)t.$$

Now observe that

$$\psi_\lambda(x) = \mathbf{E}_x[\exp(-\lambda T_D)] \geq e^{-\lambda t}\mathbf{P}_x(X(t) \in D)$$
$$\geq e^{-\lambda t}\left(1 - \sum_{i<j}\mathbf{P}_x(X_i(t) \neq X_j(t))\right).$$

This holds for any $t > 0$. Taking $t = \mathrm{dist}(x, D)$, and applying the previous estimate gives the first assertion of the lemma. For the second observe that

$$R_\lambda^0 f(x) = \mathbf{E}_x\left[\int_0^{T_D} e^{-\lambda s} f(X(s))\,ds\right] \leq \frac{\|f\|_\infty}{\lambda}(1 - \psi_\lambda(x)). \qquad \square$$

Recall that $T_\varepsilon = \inf\{t \geq 0 : |X_i(t) - X_j(t)| \geq \varepsilon \text{ for some } i,j \in \{1, 2, \ldots, N\}\}$.

LEMMA 5.2. *There exists a constant $C$, depending on $\theta$ alone, such that if $X$ is any solution to the $\mathcal{A}_N^\theta$-martingale problem starting from the diagonal, then*

$$\mathbf{E}[T_\varepsilon^2] \leq C\varepsilon^2 \qquad \text{whenever } 0 < \varepsilon < 1.$$

PROOF. Notice that $T_\varepsilon \leq \inf\{t \geq 0 : |X_1(t) - X_2(t)| \geq \varepsilon\}$. Consequently it is enough to prove the result for $N = 2$. Let $Z(t) = |X_1(t) - X_2(t)|$ and define $f(z,t) = (z^4 - 12tz^2 + 12t^2) + (2z^3 - 12tz)/\theta$ and $g(z,t) = z^2/2 - t + z/(2\theta)$ where $\theta = 2\theta(1\!:\!1)$. On applying Itô's formula we find that both $f(Z(t),t)$ and $g(Z(t),t)$ are martingales. Then using the stopping theorem we obtain first that $\mathbf{E}[T_\varepsilon] = \varepsilon^2/2 + \varepsilon/(2\theta)$, and then that

$$\mathbf{E}[T_\varepsilon^2] = 5\varepsilon^4/12 + 5\varepsilon^3/(6\theta) + \varepsilon^2/(2\theta^2)$$

from which the result follows. $\square$

LEMMA 5.3. *Under the hypothesis (5.1), the mapping $x \mapsto \psi_\lambda(x)$ is invariant under translations parallel to the diagonal, meaning $\psi_\lambda(x + y) = \psi_\lambda(x)$ for all $y \in D$. If $x \mapsto f(x)$ is invariant under such translations, so too is $x \mapsto R_\lambda^0 f(x)$.*

PROOF. Suppose $y \in D$. It is straightforward to check that if $X$ satisfies the $\mathcal{A}_N^\theta$-martingale problem starting from $x$, then $\tilde{X}(t) = X(t) + y$ solves the $\mathcal{A}_N^\theta$-martingale problem starting from $x + y$. Combining this with the uniqueness statement (5.1) we see that the stopped process $(\tilde{X}(t \wedge T_D); t \geq 0)$ under $\mathbf{P}_x$ has the same distribution as the stopped process $(X(t \wedge T_D); t \geq 0)$ under $\mathbf{P}_{x+y}$, and the lemma follows from this. $\square$

PROPOSITION 5.1.

$$\psi_\lambda(x) = \kappa_0 \mathbf{E}_x\left[\int_0^\infty e^{-\lambda s} \mathbf{1}_D(X(s))\,ds\right],$$

*where $\kappa_0$ is given by*

$$\kappa_0 = \lambda + \lim_{\varepsilon \downarrow 0} \frac{1}{\varepsilon} \sum_{v \in \mathcal{V}_+(D)} \theta(v)[1 - \psi_\lambda(\varepsilon v)].$$

*Thus under the hypothesis (5.1), for each $x \in \mathbf{R}^N$, the quantity $\mathbf{E}_x[\int_0^\infty e^{-\lambda s} \times \mathbf{1}_D(X(s))\,ds]$ does not depend on any possible choice of $\mathbf{P}_x$.*

PROOF. Introduce stopping times $T_0^\varepsilon = 0$, $T_1^\varepsilon = \inf\{t \geq 0 : X_t \in D\}$,
$$T_2^\varepsilon = \inf\{t \geq T_1^\varepsilon : |X_i(t) - X_j(t)| \geq \varepsilon \text{ for some } i,j \in \{1,2,\ldots,N\}\}$$
and in general $T_{2k+1}^\varepsilon = \inf\{t \geq T_{2k}^\varepsilon : X_t \in D\}$, and
$$T_{2k}^\varepsilon = \inf\{t \geq T_{2k-1}^\varepsilon : |X_i(t) - X_j(t)| \geq \varepsilon \text{ for some } i,j \in \{1,2,\ldots,N\}\}.$$
We have

$$
\begin{aligned}
\mathbf{E}_x&\left[\int_0^\infty e^{-\lambda s}\mathbf{1}_D(X(s))\,ds\right] \\
&= \sum_{k \text{ odd}} \mathbf{E}_x\left[\int_{T_k^\varepsilon}^{T_{k+1}^\varepsilon} e^{-\lambda s}\mathbf{1}_D(X(s))\,ds\right] \\
&= \sum_{k \text{ odd}} \mathbf{E}_x\left[e^{-\lambda T_k^\varepsilon}\tilde{\mathbf{E}}_k\left[\int_0^{T_\varepsilon} e^{-\lambda s}\mathbf{1}_D(X(s))\,ds\right]\right],
\end{aligned}
\tag{5.6}
$$

where $\tilde{\mathbf{E}}_k$ denotes expectation relative to the conditional distribution $\tilde{\mathbf{P}}_k$ of $(X(T_k^\varepsilon + u); u \geq 0)$ given $(X(u); u \leq T_k^\varepsilon)$. Notice that this holds even though it may be the case some $T_k^\varepsilon = \infty$ with positive probability, so long as we appropriately interpret the right-hand side.

Next $X$ governed by the conditional distribution $\tilde{\mathbf{P}}_k$ is almost surely also a solution to the $\mathcal{A}_N^\theta$-martingale problem, starting from the diagonal for odd $k$. Thus the estimate for $\mathbf{E}[T_\varepsilon]$ given at (3.3) is applicable, and this together with Lemmas 3.1 and 5.2 gives, for sufficiently small $\varepsilon$,

$$\left|\tilde{\mathbf{E}}_k\left[\int_0^{T_\varepsilon} e^{-\lambda s}\mathbf{1}_D(X(s))\,ds\right] - \frac{\varepsilon}{2\sum_{v \in \mathcal{V}_+(D)} \theta(v)}\right| \leq C\varepsilon^2 \quad \text{a.s.,}$$

where the constant $C$ depends on $\lambda$, $N$ and $\theta$ only. Returning to (5.6), we deduce that first, for all sufficiently small $\varepsilon$,

$$\frac{\varepsilon}{4\sum_{v \in \mathcal{V}_+(D)} \theta(v)} \sum_{k \text{ odd}} \mathbf{E}_x[e^{-\lambda T_k^\varepsilon}] \leq \mathbf{E}_x\left[\int_0^\infty e^{-\lambda s}\mathbf{1}_D(X(s))\,ds\right] < \infty$$

and then that

$$\mathbf{E}_x\left[\int_0^\infty e^{-\lambda s}\mathbf{1}_D(X(s))\,ds\right] = \lim_{\varepsilon \downarrow 0} \frac{\varepsilon}{2\sum_{v \in \mathcal{V}_+(D)} \theta(v)} \sum_{k \text{ odd}} \mathbf{E}_x[e^{-\lambda T_k^\varepsilon}]. \tag{5.7}$$

Next we use the same sequence of stopping times $T_k^\varepsilon$ to decompose $\psi_\lambda(x)$:

$$\psi_\lambda(x) = \lim_{n \to \infty} \mathbf{E}_x[\psi_\lambda(x) - e^{-\lambda T_n^\varepsilon}\psi_\lambda(X(T_n^\varepsilon))]$$

$$\text{(5.8)} \quad = \lim_{n \to \infty} \mathbf{E}_x\left[\sum_{k=0}^n e^{-\lambda T_k^\varepsilon}\psi_\lambda(X(T_k^\varepsilon)) - e^{-\lambda T_{k+1}^\varepsilon}\psi_\lambda(X(T_{k+1}^\varepsilon))\right]$$

$$= \sum_{k \text{ odd}} \mathbf{E}_x[e^{-\lambda T_k^\varepsilon}\tilde{\mathbf{E}}_k[1 - e^{-\lambda T_\varepsilon}\psi_\lambda(X(T_\varepsilon))]].$$

Here the even terms of the sum are dropped because $\tilde{\mathbf{E}}_k[\psi_\lambda(X(0)) - \exp(-\lambda \times T_D)] = 0$ for each even $k$, by virtue of the uniqueness hypothesis (5.1) and the definition of $\psi_\lambda$.

We claim that there exists a constant $C$ depending on $\lambda$, $N$ and $\theta$ so that

$$\text{(5.9)} \quad \left|\tilde{\mathbf{E}}_k[1 - e^{-\lambda T_\varepsilon}\psi_\lambda(X(T_\varepsilon))] - \frac{\lambda \varepsilon/2 + \sum_{v \in \mathcal{V}_+(D)}[1 - \psi_\lambda(\varepsilon v/2)]\theta(v)}{\sum_{u \in \mathcal{V}_+(D)}\theta(u)}\right|$$
$$\leq C\varepsilon^{3/2} \qquad \text{a.s.}$$

Taken together with (5.7) and (5.8) this implies that the limit defining $\kappa_0$ exists, and that the equality in the statement of the proposition holds.

In verifying the claimed inequality we will write $\mathbf{E}$ instead of $\tilde{\mathbf{E}}_k$, our arguments holding for any solution to the $\mathcal{A}_N^\theta$-martingale problem starting from the diagonal. In the following $C_1, C_2, \ldots$ denote constants that may depend on $\theta$, $N$ and $\lambda$ only. We begin with

$$\text{(5.10)} \quad \mathbf{E}[1 - \psi_\lambda(X(T_\varepsilon))] = \sum_{v \in \mathcal{V}_+(D)}(1 - \psi_\lambda(\varepsilon v/2))\mathbf{P}(X(T_\varepsilon) \in E(v))$$
$$+ \mathbf{E}[(1 - \psi_\lambda(X(T_\varepsilon)))\mathbf{1}_{\Lambda(\varepsilon)}],$$

where $\Lambda(\varepsilon)$ denotes the event that at least three of the coordinates of $X(T_\varepsilon)$ take distinct values. Here we have used the fact that $\psi_\lambda(X(T_\varepsilon)) = \psi_\lambda(\varepsilon v/2)$ on the event $X(T_\varepsilon) \in E(v)$ which follows from Lemma 5.3. Now from Proposition 3.1 we know that $\mathbf{P}(\Lambda(\varepsilon)) \leq C_1\varepsilon$, and from Lemma 5.1 that $(1 - \psi_\lambda(X(T_\varepsilon))) \leq C_2\varepsilon^{1/2}$, and thus $\mathbf{E}[(1 - \psi_\lambda(X(T_\varepsilon)))\mathbf{1}_\Lambda] \leq C_3\varepsilon^{3/2}$. Applying the estimate (3.4) for $\mathbf{P}(X(T_\varepsilon) \in E(v))$ with another appeal to Lemma 5.1, we deduce from (5.10) that

$$\text{(5.11)} \quad \left|\mathbf{E}[1 - \psi_\lambda(X(T_\varepsilon))] - \sum_{v \in \mathcal{V}_+(D)}[1 - \psi_\lambda(\varepsilon v/2)]\frac{\theta(v)}{\sum_{u \in \mathcal{V}_+(D)}\theta(u)}\right| \leq C_5\varepsilon^{3/2}.$$

Next we have

$$\text{(5.12)} \quad |\mathbf{E}[1 - e^{-\lambda T_\varepsilon}\psi_\lambda(X(T_\varepsilon))] - \mathbf{E}[1 + \lambda T_\varepsilon - \psi_\lambda(X(T_\varepsilon))]|$$
$$\leq \lambda \mathbf{E}[T_\varepsilon(1 - \psi_\lambda(X(T_\varepsilon)))] + \lambda^2 \mathbf{E}[T_\varepsilon^2 \psi_\lambda(X(T_\varepsilon))] \leq C_6\varepsilon^{3/2}$$

using Lemmas 5.1 and 5.2 together with the estimate for $\mathbf{E}[T_\varepsilon]$ given at (3.3). Finally the desired inequality (5.9) follows from putting together (5.11), (5.12) and using (3.3) once again. $\square$

PROPOSITION 5.2. *Suppose that $f$ is bounded, zero in a neighborhood of $D$, and invariant under shifts along $D$: that is, $f(x+y) = f(x)$ for all $y \in D$. Then*

$$\mathbf{E}_x\left[\int_0^\infty e^{-\lambda s} f(X(s))\,ds\right] = R_\lambda^0 f(x) + \frac{\kappa_f}{\kappa_0}\psi_\lambda(x),$$

*where*

$$\kappa_f = \lim_{\varepsilon\downarrow 0}\frac{1}{\varepsilon}\sum_{v\in\mathcal{V}_+(D)}\theta(v)R_\lambda^0 f(\varepsilon v).$$

*Thus, under the hypothesis (5.1), for each $x \in \mathbf{R}^N$, the quantity $\mathbf{E}_x[\int_0^\infty e^{-\lambda s}\times f(X(s))\,ds]$ does not depend on any possible choice of $\mathbf{P}_x$.*

PROOF. Suppose that $\varepsilon$ is small enough that $f(x) = 0$ for all $x \in \mathbf{R}^N$ within a distance $N\varepsilon$ of $D$. Then using the same sequence of stopping times as in the previous proof,

$$R_\lambda^0 f(x) = \lim_{n\to\infty}\mathbf{E}_x[R_\lambda^0 f(x) - e^{-\lambda T_n^\varepsilon}R_\lambda^0 f(X(T_n^\varepsilon))]$$

$$= \lim_{n\to\infty}\mathbf{E}_x\left[\sum_{k=0}^n e^{-\lambda T_k^\varepsilon}R_\lambda^0 f(X(T_k^\varepsilon)) - e^{-\lambda T_{k+1}^\varepsilon}R_\lambda^0 f(X(T_{k+1}^\varepsilon))\right]$$

$$= \sum_{k\text{ even}}\mathbf{E}_x\left[\int_{T_k^\varepsilon}^{T_{k+1}^\varepsilon}e^{-\lambda s}f(X(s))\,ds\right] - \sum_{k\text{ odd}}\mathbf{E}_x[e^{-\lambda T_{k+1}^\varepsilon}R_\lambda^0 f(X(T_{k+1}^\varepsilon))]$$

$$= \mathbf{E}_x\left[\int_0^\infty e^{-\lambda s}f(X(s))\,ds\right] - \sum_{k\text{ odd}}\mathbf{E}_x[e^{-\lambda T_k^\varepsilon}\tilde{\mathbf{E}}_k[e^{-\lambda T_\varepsilon}R_\lambda^0 f(X(T_\varepsilon))]].$$

Here we have made use of the equality that for even $k$,

$$\mathbf{E}_x\left[\int_{T_k^\varepsilon}^{T_{k+1}^\varepsilon}e^{-\lambda s}f(X(s))\,ds\right] = \mathbf{E}_x\left[e^{-\lambda T_k^\varepsilon}\tilde{\mathbf{E}}_k\left[\int_0^{T_D}e^{-\lambda s}f(X(s))\,ds\right]\right]$$

$$= \mathbf{E}_x[e^{-\lambda T_k^\varepsilon}R_\lambda^0 f(X(T_k^\varepsilon))],$$

which holds since $X$ governed by the conditional distribution $\tilde{\mathbf{P}}_k$ is also a solution to the $\mathcal{A}_N^\theta$-martingale problem.

We claim that there exists a constant $C$ not depending on $k$ so that

$$\left|\tilde{\mathbf{E}}_k[e^{-\lambda T_\varepsilon}R_\lambda^0 f(X(T_\varepsilon))] - \sum_{v\in\mathcal{V}_+(D)}R_\lambda^0 f(\varepsilon v/2)\frac{\theta(v)}{\sum_{u\in\mathcal{V}_+(D)}\theta(u)}\right| \leq C\varepsilon^{3/2} \quad \text{a.s.}$$

Taken together with the preceding decomposition of $R_\lambda^0 f(x)$ and (5.7), this implies that the limit defining $\kappa_f$ exists, and that the equality in the statement of the proposition holds.

As in the previous proposition we write $\mathbf{E}$ instead of $\tilde{\mathbf{E}}_k$, our arguments holding for any solution to the $\mathcal{A}_N^\theta$-martingale problem starting from the diagonal. $C_1, C_2, \ldots$ denote constants that may depend on $\theta$, $N$ and $\lambda$ only. We begin with

$$\mathbf{E}[R_\lambda^0 f(X(T_\varepsilon))] = \sum_{v \in \mathcal{V}_+(D)} R_\lambda^0 f(\varepsilon v/2) \mathbf{P}(X(T_\varepsilon) \in E(v))$$

(5.13)
$$+ \mathbf{E}[R_\lambda^0 f(X(T_\varepsilon)) \mathbf{1}_{\Lambda(\varepsilon)}].$$

Here we have used the fact that $R_\lambda^0 f(X(T_\varepsilon)) = R_\lambda^0 f(\varepsilon v/2)$ on the event $X(T_\varepsilon) \in E(v)$ which follows from Lemma 5.3.

Combining Proposition 3.1 and Lemma 5.1 we obtain $\mathbf{E}[R_\lambda^0 f(X(T\varepsilon))\mathbf{1}_\Lambda] \leq C_1 \varepsilon^{3/2}$. Applying the estimate (3.4) for $\mathbf{P}(X(T_\varepsilon) \in E(v))$, with another appeal to Lemma 5.1, we deduce from (5.13) that

(5.14) $$\left| \mathbf{E}[R_\lambda^0 f(X(T_\varepsilon))] - \sum_{v \in \mathcal{V}_+(D)} R_\lambda^0 f(\varepsilon v/2) \frac{\theta(v)}{\sum_{u \in \mathcal{V}_+(D)} \theta(u)} \right| \leq C_2 \varepsilon^{3/2}.$$

To finish we have

$$|\mathbf{E}[e^{-\lambda T_\varepsilon} R_\lambda^0 f(X(T_\varepsilon))] - \mathbf{E}[R_\lambda^0 f(X(T_\varepsilon))]| \leq C_3 \varepsilon^{3/2},$$

by combining Lemma 5.1 with an upper bound for $\mathbf{E}[T_\varepsilon]$ of order $\varepsilon$ which follows from (3.3). $\square$

For a point $x \in \mathbf{R}^N$ we denote by $\hat{x}$ the point obtained by orthogonal projection of $x$ onto the hyperplane $\{x \in \mathbf{R}^N : \sum x_i = 0\}$.

PROPOSITION 5.3. *Suppose the hypothesis (5.1) holds. Then if $X$ is any solution to the $\mathcal{A}_N^\theta$-martingale problem starting from a given point $x \in \mathbf{R}^N$, the law of the projected process $(\hat{X}(t); t \geq 0)$ is uniquely determined.*

PROOF. Suppose that $\tilde{\mathbf{P}}_x$ and $\mathbf{P}_x$ are two probability measures under each of which $X$ is a solution to the $\mathcal{A}_N^\theta$-martingale problem starting from $x$. Then by the previous proposition if $f$ is invariant under shifts along the diagonal $D$, and zero in a neighborhood of $D$,

$$\tilde{\mathbf{E}}_x \left[ \int_0^\infty e^{-\lambda s} f(X(s)) \, ds \right] = \mathbf{E}_x \left[ \int_0^\infty e^{-\lambda s} f(X(s)) \, ds \right].$$

But both sides are positive measures (integrating $f$) with the same total mass, so by a monotone convergence argument, this equality holds without the assumption of $f$ being zero in a neighborhood of $D$.

Next by inverting the Laplace transform, and using continuity in $t$, we deduce that

$$\tilde{\mathbf{E}}_x[f(X(t))] = \mathbf{E}_x[f(X(t))] \qquad \text{for all } t \geq 0.$$

Notice that $x \mapsto \mathbf{E}_x[f(X(t))]$ is itself invariant under shifts along $D$. Using this, and the fact that the conditional distribution of process $(X(t+u); u \geq 0)$ given $(X(s); s \leq t)$ solves the martingale problem, we obtain in a standard way

$$\tilde{\mathbf{E}}_x\left[\prod_{i=1}^n f_i(X(t_i))\right] = \mathbf{E}_x\left[\prod_{i=1}^n f_i(X(t_i))\right]$$

for any bounded functions $f_1, f_2, \ldots, f_n$, and times $0 < t_1 < t_2 < \cdots < t_n$, with each $f_k$ invariant under shifts. The claimed result follows from this. □

In order to extend the uniqueness result from $\hat{X}$ to $X$ we use a change of measure technique. We now assume that $\theta(0:1) = \theta(1:0) = 0$ and consequently $X$ is a martingale. For $x \in \mathbf{R}^N$ and $i = 1, 2, \ldots, N$ let $m_i(x)$ be the cardinality of the set $\{k \in \{1, 2, \ldots, n\} : x_k = x_i\}$. Then we define a martingale $(M(t); t \geq 0)$ via

(5.15) $$M(t) = \sum_{i=1}^N \int_0^t \frac{1}{m_i(X(s))} dX_i(s).$$

It is easy to check that $\langle M, X_i \rangle(t) = t$ for each $i$, and consequently $M$ and $\hat{X}$ are orthogonal as martingales. A similar calculation also gives $\langle M \rangle(t) = \int_0^t |\pi(X(s))| \, ds$, where $|\pi(x)|$ is the number of components in the partition $\pi(x)$ associated with $x$.

Suppose $(\beta(t); t \geq 0)$ is a bounded measurable, nonrandom process, and define $Z$ to be the exponential local martingale

(5.16) $$Z(t) = \exp\left\{\int_0^t \beta(s) \, dM(s) - \frac{1}{2}\int_0^t \beta(s)^2 d\langle M \rangle(s)\right\}.$$

By Novikov's criterion this is a true martingale. In the following **1** denotes the vector $(1, 1, \ldots, 1) \in \mathbf{R}^N$.

LEMMA 5.4. *Suppose $(X(t); t \geq 0)$, governed by a probability measure* **P**, *solves the $\mathcal{A}_N^\theta$-martingale problem. Let $\tilde{\mathbf{P}}$ be locally absolutely continuous with respect to* **P** *with density*

$$\tilde{\mathbf{P}} = Z(t) \cdot \mathbf{P} \qquad \text{on } \mathcal{F}_t,$$

*where $(\mathcal{F}_t; t \geq 0)$ is the filtration generated by $X$. Then under $\tilde{\mathbf{P}}$,*

$$X(t) - \int_0^t \beta(s)\mathbf{1} \, ds \qquad \text{solves the } \mathcal{A}_N^\theta\text{-martingale problem.}$$

PROOF. Let $\tilde{X}(t) = X(t) - \int_0^t \beta(s)\mathbf{1}\,ds$. Evidently $\tilde{X}$ has the correct quadratic variation. We must show that for each $f \in L_N$ the process

$$\tilde{M}^f(t) = f(\tilde{X}(t)) - \int_0^t \mathcal{A}_N^\theta f(\tilde{X}(s))\,ds$$

is a martingale under $\tilde{P}$.

As in the proof of Proposition 2.2 we note that $\tilde{X}(t)$ and $X(t)$ always belong to the same cell, and by the linearity of $f$, and the constancy of $\mathcal{A}_N^\theta f$ on each cell,

$$\tilde{M}^f(t) = f(X(t)) - \int_0^t \beta(s)f(\mathbf{1})\,ds - \int_0^t \mathcal{A}_N^\theta f(X(s))\,ds.$$

$M^f(t)$ can be written as a stochastic integral with respect to $X$, and using this and the observation that $\langle X_i, M\rangle(t) = t$ made previously we are able to calculate that $\langle M^f, M\rangle(t) = f(\mathbf{1})t$. Consequently,

$$\langle M^f, M^\beta\rangle(t) = \int_0^t \beta(s)f(\mathbf{1})\,ds,$$

where $M^\beta(t)$ denotes the martingale $\int_0^t \beta(s)\,dM(s)$. By Girsanov's theorem the process

$$M^f(t) - \langle M^f, M^\beta\rangle(t),$$

is a martingale under $\tilde{P}$, and so we are done. $\square$

PROPOSITION 5.4. *Suppose the hypothesis (5.1) holds. Then if $X$ is any solution to the $\mathcal{A}_N^\theta$-martingale problem starting from a given point $x \in \mathbf{R}^N$, the law of the process $(X(t); t \geq 0)$ is uniquely determined.*

PROOF. Notice that the process $X$ can be expressed as some measurable function of the pair of processes $(\hat{X}, M)$ since

$$\frac{1}{N}\sum_{i=1}^N X_i(t) = \int_0^t \left\{\sum m_j(\hat{X}(s))^{-1}\right\}^{-1}\left\{dM(s) - \sum_i \frac{d\hat{X}_i(s)}{m_i(\hat{X}(s))}\right\}.$$

Thus, since we know the law of $\hat{X}$ is uniquely determined, it is enough to show that the conditional distribution of $M$ given $\hat{X}$ is unique also. For this we argue as follows, using the notation introduced in the previous lemma. Let $T > 0$. For any bounded function $F: \mathbf{C}([0,T], \mathbf{R}^N) \to \mathbf{R}$, we have by definition of $\tilde{\mathbf{P}}$, $\tilde{\mathbf{E}}[F(\hat{X})] = \mathbf{E}[Z(T)F(\hat{X})]$. By the previous lemma, the process $\tilde{X}$, under $\tilde{\mathbf{P}}$, solves the $\mathcal{A}_N^\theta$-martingale problem. So applying the uniqueness result Proposition 5.3 we deduce that $\tilde{\mathbf{E}}[F(\hat{X})] = \mathbf{E}[F(\hat{X})]$. From the resulting equality, $\mathbf{E}[Z(T)F(\hat{X})] = \mathbf{E}[F(\hat{X})]$, $F$ being arbitrary, it follows that

$$\mathbf{E}\left[\exp\left\{\int_0^T \beta(s)\,dM(s)\right\}\Big|(\hat{X}(t); t\in[0,T])\right] = \exp\left\{\frac{1}{2}\int_0^T \beta(s)^2|\pi(\hat{X}(s))|\,ds\right\}.$$

The process $\beta$ and time $T$ being arbitrary, this shows that the distribution of $M$ given $\hat{X}$ is uniquely specified. □

Finally the proof of the uniqueness statement of Theorem 2.1 is completed by an argument by induction on the dimension $N$. We can assume that $\theta(0:1) = \theta(1:0) = 0$. Uniqueness for $N=1$ holds by Lévy's characterization of Brownian motion. Then assuming uniqueness holds for every dimension $n \leq N-1$, Proposition 4.2 asserts that the hypothesis (5.1) is met, and thus by Proposition 5.4, uniqueness also holds for dimension $N$.

**6. An approximation scheme and the proof of existence.** In this section we construct a solution to the $\mathcal{A}_N^\theta$-martingale problem, as the scaling limit of a sequence of Markov chains.

Let $p = (p(k:l); k,l \geq 0)$ be a family of nonnegative parameters satisfying the consistency condition

(6.1) $\qquad p(k:l) = p(k+1:l) + p(k:l+1) \qquad$ for all $k,l \geq 0$.

We will consider a continuous time Markov chain $(Y(t); t \geq 0)$ with state space the integer lattice $\mathbf{Z}^N$ which has generator given by

(6.2) $\qquad\qquad \mathcal{G}_N^p f(x) = \sum_{v \in \mathcal{V}(x)} p(v)\{f(x+v) - f(x)\}.$

Here we write $p(v)$ for $p(|I|,|J|)$, where $v = v_{IJ}$.

The following proposition is proved exactly as Proposition 2.1.

PROPOSITION 6.1. *Suppose that $Y$ is a Markov chain with generator $\mathcal{G}_N^p$, and let $Z$ be the process consisting of some $M$ coordinates of $Y$. Then $Z$ is a Markov chain with generator $\mathcal{G}_M^p$.*

From now on we assume that $p$ satisfies $p(0:0) = p(1:0) + p(0:1) = 1$. Define $d = p(1:0) - p(0:1)$. Then with the help of the preceding proposition we see that each coordinate $Y_i$ of $Y$ is a simple random walk on $\mathbf{Z}$ with drift $d$, and in particular

(6.3) $\qquad Y_i(t) - dt \quad$ and $\quad (Y_i(t) - dt)^2 - t \qquad$ are martingales

relative to the natural filtration of $Y$. Similarly, we can consider any pair of coordinates $(Y_i, Y_j)$ and we find that they evolve independently from each other when apart but interact when they meet. In fact

(6.4) $\quad (Y_i(t) - dt)(Y_j(t) - dt) - (1 - 4p(1\!:\!1)) \int_0^t \mathbf{1}(Y_i(s) = Y_j(s))\, ds \quad$ and

(6.5) $\quad |Y_i(t) - Y_j(t)| - 4p(1\!:\!1) \int_0^t \mathbf{1}(Y_i(s) = Y_j(s))\, ds$

are both martingales.

Fix the integer $N \geq 1$. Let $(p_n; n \geq 1)$ be a sequence of families of parameters, all satisfying (6.1) and $p_n(0\!:\!0) = 1$, and such that as $n$ tends to infinity,

(6.6) $\qquad n^{1/2}(p_n(k\!:\!l) - \tfrac{1}{2}\mathbf{1}(k=0) - \tfrac{1}{2}\mathbf{1}(l=0)) \to \theta(k\!:\!l)$

for all $0 \leq k, l \leq N$, where $(\theta(k\!:\!l); k, l \geq 0)$ satisfies the consistency and positivity conditions (2.7) and (2.9). Let $(x_n; n \geq 0)$ be a sequence of points in $\mathbf{R}^N$ converging to a point $x$, with $x_n \in n^{-1/2}\mathbf{Z}^N$ for every $n$. For $n \geq 1$, let $Y^n$ be the scaled process given by $Y^n(t) = n^{-1/2}Y(nt)$ for $t \geq 0$, where $Y$ is a Markov chain with generator $\mathcal{G}_N^{p_n}$ starting from $n^{1/2}x_n$.

In the following convergence in law means weak convergence of probability measures on the Skorokhod space $\mathbf{D}([0, \infty), \mathbf{R}^N)$.

PROPOSITION 6.2. *Suppose that the sequence of processes $(Y^n(t); t \geq 0)$ converges in law to a process $(X(t); t \geq 0)$. Then for each $i \neq j$ the processes $(X_i(t) - 2\beta t; t \geq 0)$ and $(X_j(t) - 2\beta t; t \geq 0)$ are a pair of $\theta$-coupled Brownian motions where $\theta = 2\theta(1\!:\!1)$ and $\beta = \theta(1\!:\!0) = -\theta(0\!:\!1)$ are determined by (6.6).*

PROOF. Notice that the drift of $Y^n$ is given by $d_n = n^{1/2}(p_n(1\!:\!0) - p_n(0\!:\!1))$ which tends to $2\beta$ as $n$ tends to infinity. Thus by the usual central limit theorem for random walks each $X_i(t) - 2\beta t$ is a Brownian motion starting from $x_i$. Moreover it is a martingale relative to the natural filtration of $X$, since the corresponding statement holds for the coordinates of each $Y^n$.

Our first task is to determine the quadratic covariation of $X_i$ and $X_j$. Consider two times $0 \leq t_1 \leq t_2$. Let $g : \mathbf{D}([0, \infty), \mathbf{R}^N) \to \mathbf{R}$ be nonnegative, bounded, continuous and measurable with respect to $\mathcal{D}_{t_1}$, where $(\mathcal{D}_t; t \geq 0)$ is the filtration generated by the coordinate process. The mapping $\alpha \mapsto \int_{t_1}^{t_2} \mathbf{1}(\alpha_i(s) = \alpha_j(s))\, ds$ is upper semicontinuous relative to the Skorohod topology on $\mathbf{D}([0, \infty), \mathbf{R}^N)$. Thus by weak convergence,

(6.7)
$$\mathbf{E}\!\left[g(X) \int_{t_1}^{t_2} \mathbf{1}(X_i(s) = X_j(s))\, ds\right]$$
$$\geq \limsup_{n \to \infty} \mathbf{E}\!\left[g(Y^n) \int_{t_1}^{t_2} \mathbf{1}(Y_i^n(s) = Y_j^n(s))\, ds\right].$$

Let us put $\hat{Y}^n(t) = Y^n(t) - d_n t$ and $\hat{X}(t) = X(t) - 2\beta t$. It is easily verified that the family of random variables $(\hat{Y}_i^n(t_2)\hat{Y}_j^n(t_2) - \hat{Y}_i^n(t_1)\hat{Y}_j^n(t_1); n \geq 1)$ is uniformly integrable, and weak convergence gives

$$
\mathbf{E}[g(X)(\hat{X}_i(t_2)\hat{X}_j(t_2) - \hat{X}_i(t_1)\hat{X}_j(t_1))]
$$
(6.8)
$$
= \lim_{n \to \infty} \mathbf{E}[g(Y^n)(\hat{Y}_i^n(t_2)\hat{Y}_j^n(t_2) - \hat{Y}_i^n(t_1)\hat{Y}_j^n(t_1))].
$$

Now from (6.4), account taken of the scaling and for $n$ large enough that $p_n(1\!:\!1) \leq 1/4$, we deduce that

$$
0 \leq \mathbf{E}[g(Y^n)(\hat{Y}_i^n(t_2)\hat{Y}_j^n(t_2) - \hat{Y}_i^n(t_1)\hat{Y}_j^n(t_1))]
$$
$$
\leq \mathbf{E}\bigg[g(Y^n) \int_{t_1}^{t_2} \mathbf{1}(Y_i^n(s) = Y_j^n(s))\,ds\bigg].
$$

Combining this with the above consequences of weak convergence we obtain

$$
0 \leq \mathbf{E}[g(X)(\hat{X}_i(t_2)\hat{X}_j(t_2) - \hat{X}_i(t_1)\hat{X}_j(t_1))] \leq \mathbf{E}\bigg[g(X) \int_{t_1}^{t_2} \mathbf{1}(X_i(s) = X_j(s))\,ds\bigg].
$$

This being true for arbitrary $t_1 \leq t_2$ and $g$, it follows that the bracket $\langle X_i, X_j \rangle$ must be an increasing process such that $\langle X_i, X_j \rangle(t) - \int_0^t \mathbf{1}(X_i(s) = X_j(s))\,ds$ is a decreasing process, in particular the measure $d\langle X_i, X_j \rangle(t)$ is carried by the set $\{t : X_i(t) = X_j(t)\}$. But the occupation time formula for semimartingales, see Chapter VI of [13], implies that $d\langle X_i - X_j \rangle(t)$ does not charge the set $\{t : X_i(t) - X_j(t) = 0\}$ and using the fact that $\langle X_i \rangle(t) + \langle X_j \rangle(t) = 2t$ we are able to conclude that $\langle X_i, X_j \rangle(t) = \int_0^t \mathbf{1}(X_i(s) = X_j(s))\,ds$. Using this and (6.8) we can strengthen (6.7) to

$$
\mathbf{E}\bigg[g(X) \int_{t_1}^{t_2} \mathbf{1}(X_i(s) = X_j(s))\,ds\bigg]
$$
(6.9)
$$
= \lim_{n \to \infty} \mathbf{E}\bigg[g(Y^n) \int_{t_1}^{t_2} \mathbf{1}(Y_i^n(s) = Y_j^n(s))\,ds\bigg].
$$

To prove the proposition it remains to show that the process

$$
|X_i(t) - X_j(t)| - 2\theta \int_0^t \mathbf{1}(X_i(s) = X_j(s))\,ds
$$

is a martingale. With $t_1 \leq t_2$ and $g$ as before, by weak convergence and uniform integrability,

$$
\mathbf{E}[g(X)(|X_i(t_2) - X_j(t_2)| - |X_i(t_1) - X_j(t_1)|)]
$$
(6.10)
$$
= \lim_{n \to \infty} \mathbf{E}[g(Y^n)(|Y_i^n(t_2) - Y_j^n(t_2)| - |Y_i^n(t_1) - Y_j^n(t_1)|)].
$$

From (6.5) we know that

$$\mathbf{E}[g(Y^n)(|Y_i^n(t_2) - Y_j^n(t_2)| - |Y_i^n(t_1) - Y_j^n(t_1)|)]$$
$$= 4n^{1/2} p_n(1:1) \mathbf{E}\left[g(Y^n) \int_{t_1}^{t_2} \mathbf{1}(Y_i^n(s) = Y_j^n(s)) \, ds\right].$$

Combining this with (6.9), (6.10) and the fact $4n^{1/2} p_n(1:1)$ tends to $2\theta$ as $n$ tends to infinity, we deduce

$$\mathbf{E}[g(X)(|X_i(t_2) - X_j(t_2)| - |X_i(t_1) - X_j(t_1)|)]$$
$$= 2\theta \mathbf{E}\left[g(X) \int_{t_1}^{t_2} \mathbf{1}(X_i(s) = X_j(s)) \, ds\right],$$

which completes the proof. □

PROPOSITION 6.3. *Suppose that the sequence of processes $(Y^n(t); t \geq 0)$ converges in law to a process $(X(t); t \geq 0)$. Then the process $X$ solves the $\mathcal{A}_N^\theta$-martingale problem, where $\theta$ is determined by (6.6).*

PROOF. Suppose that $f \in L_N$ and $x \mapsto \mathcal{A}_N^\theta f(x)$ is upper semicontinuous on $\mathbf{R}^N$. We will show that, in this case, the process $M_t^f = f(X(t)) - \int_0^t \mathcal{A}_N^\theta f(X(s)) \, ds$ is a supermartingale.

We begin by observing that we may re-write the generator $\mathcal{G}_N^{p_n}$ in a form closer to that of $\mathcal{A}_N^\theta$. In fact we have, for $x \in \mathbf{Z}^N$,

$$\mathcal{G}_N^{p_n} f(x) = \sum_{v \in \mathcal{V}(x)} p_n(v) \{f(x+v) - f(x)\} = \sum_{v \in \mathcal{V}(x)} p_n(v) \nabla_v f(x).$$

This is because $f$ is linear when restricted to any cone of the form $K_\rho = \{x \in \mathbf{R}^N : x_i \leq x_j \text{ if } \rho(i) < \rho(j)\}$ for $\rho$ a permutation on $\{1, 2, \ldots, N\}$, and for $x \in \mathbf{Z}^N$ the points $x$ and $x+v$ always belong to some common cone. Also, since $f$ is linear on each cell, and thus $\nabla_v f(x) = -\nabla_{-v} f(x)$ for any $v \in \mathcal{V}_0(x)$, we see that, for $x \in \mathbf{Z}^N$,

$$\mathcal{G}_N^{p_n} f(x) = \sum_{v \in \mathcal{V}(x)} \hat{p}_n(v) \nabla_v f(x),$$

where $\hat{p}(k:l) = p(k:l) - \frac{1}{2}\mathbf{1}(k=0) - \frac{1}{2}\mathbf{1}(l=0)$. From this, using the convergence of $p_n$ given by (6.6), and noting that $\mathcal{A}_N^\theta f(x/\sqrt{n}) = \mathcal{A}_N^\theta f(x)$, we obtain

$$(6.11) \qquad \sup_{x \in \mathbf{Z}^N} |n^{1/2} \mathcal{G}_N^{p_n} f(x) - \mathcal{A}_N^\theta f(x/\sqrt{n})| \to 0$$

as $n$ tends to infinity.

Consider two times $0 \leq t_1 \leq t_2$. Let $g : \mathbf{D}([0, \infty), \mathbf{R}^N) \to \mathbf{R}$ be nonnegative, bounded, continuous and measurable with respect to $\mathcal{D}_{t_1}$. Now the mapping $\alpha \mapsto \int_{t_1}^{t_2} \mathcal{A}_N^\theta f(\alpha(s))\,ds$ is upper semicontinuous relative to the Skorohod topology on $\mathbf{D}([0,\infty),\mathbf{R}^N)$, as a consequence of $x \mapsto \mathcal{A}_N^\theta f(x)$ being upper semicontinuous on $\mathbf{R}^N$. It is also, for fixed $t_1$ and $t_2$, bounded. Thus by weak convergence, $g$ being bounded, continuous and nonnegative,

$$\mathbf{E}\left[g(X)\int_{t_1}^{t_2} \mathcal{A}_N^\theta f(X(s))\,ds\right] \geq \limsup_{n\to\infty} \mathbf{E}\left[g(Y^n)\int_{t_1}^{t_2} \mathcal{A}_N^\theta f(Y^n(s))\,ds\right]$$
$$= \limsup_{n\to\infty} n^{1/2}\mathbf{E}\left[g(Y^n)\int_{t_1}^{t_2} \mathcal{G}_N^{p_n} f(\sqrt{n}Y^n(s))\,ds\right],$$

where the equality is due to the uniform convergence (6.11). The family of random variables $(f(Y^n(t_2)) - f(Y^n(t_1)); n \geq 1)$ is uniformly integrable because

$$\mathbf{E}[(f(Y^n(t_2)) - f(Y^n(t_1)))^2] \leq \|f\|_{\mathrm{Lip}}^2 \sum_{i=1}^N \mathbf{E}[(Y_i^n(t_2) - Y_i^n(t_1))^2]$$
$$= N\|f\|_{\mathrm{Lip}}^2(d_n^2(t_2 - t_1)^2 + (t_2 - t_1)),$$

where $\|f\|_{\mathrm{Lip}}$ is the Lipshitz norm of $f$ and $d_n = n^{1/2}(p_n(1\!:\!0) - p_n(0\!:\!1)) \to 2\beta$ as $n$ tends to infinity. Thus by weak convergence

$$\mathbf{E}[g(X)(f(X(t_2)) - f(X(t_1)))] = \lim_{n\to\infty} \mathbf{E}[g(Y^n)(f(Y^n(t_2)) - f(Y^n(t_1)))].$$

Next, account taken of the scaling, we find that

$$M_t^{f,n} = f(Y^n(t)) - n^{1/2}\int_0^t \mathcal{G}_N^{p_n} f(\sqrt{n}Y^n(s))\,ds \quad \text{is a martingale}$$

relative to the filtration of $Y^n$, and so,

$$\mathbf{E}[g(Y^n)(f(Y^n(t_2)) - f(Y^n(t_1)))] = n^{1/2}\mathbf{E}\left[g(Y^n)\int_{t_1}^{t_2} \mathcal{G}_N^{p_n} f(\sqrt{n}Y^n(s))\,ds\right].$$

Combining this with the consequences of weak convergence obtained above we deduce that

$$\mathbf{E}[g(X)(f(X(t_2)) - f(X(t_1)))] \leq \mathbf{E}\left[g(X)\int_{t_1}^{t_2} \mathcal{A}_N^\theta f(X(s))\,ds\right].$$

This proves that $M_t^f$ is a supermartingale as claimed.

Now consider a general $f \in L_N$. Let $g : \mathbf{R}^N \to \mathbf{R}$ be given by $g(x) = \sum_{i \neq j} |x_i - x_j|$. Then for sufficiently large $R > 0$, both $\mathcal{A}_N^\theta(Rg + f)$ and $\mathcal{A}_N^\theta(Rg - f)$ are upper semicontinuous. Applying the above to $Rg + f$ and $Rg - f$ in turn we deduce that

(6.12)
$$(Rg + f)(X(t)) - \int_0^t \mathcal{A}_N^\theta(Rg + f)(X(s))\,ds \quad \text{and}$$
$$(Rg - f)(X(t)) - \int_0^t \mathcal{A}_N^\theta(Rg - f)(X(s))\,ds$$

are supermartingales. But

$$g(X(t)) - \int_0^t \mathcal{A}_N^\theta g(X(s)) \, ds$$
$$= \sum_{i \neq j} \left\{ |X_i(t) - X_j(t)| - 4\theta(1\!:\!1) \int_0^t \mathbf{1}_{(X_i(s)=X_j(s))} \, ds \right\}$$

is a martingale by the preceding proposition. Subtracting $R$ times this martingale from the two supermartingales at (6.12) we conclude that $f(X(t)) - \int_0^t \mathcal{A}_N^\theta f(X(s)) \, ds$ is a martingale. □

Of course in the above proposition we may replace the sequence of processes $(Y^n; n \geq 1)$ by any subsequence $(Y^{n_r}; r \geq 1)$ and the conclusion still holds. Thus to prove the existence to a solution of the $\mathcal{A}_N^\theta$-martingale problem, we must observe that $(Y^n; n \geq 1)$ is tight, and consequently a subsequence converging in law does exist. But this tightness follows from the fact that for each $i$, the sequence of processes $(Y_i^n; n \geq 1)$ is converging in law to a process with continuous paths; see VI.3.33 of [8].

Finally the proof of the existence statement of Theorem 2.1 is complete once we have observed that every family of parameters $\theta$ can be obtained from some suitable sequence of parameters $(p_n; n \geq 1)$. This may be achieved, for a given $\theta$, by setting

(6.13) $\quad p_n(k\!:\!l) = \frac{1}{2}\mathbf{1}(k=0) + \frac{1}{2}\mathbf{1}(l=0) + n^{-1/2}\theta(k\!:\!l),$

noting that for a fixed $N$ and all sufficiently large $n$, this gives $p_n(k\!:\!l) \geq 0$, at least for $0 \leq k, l \leq N$, which is sufficient.

**7. On a diffusion in the half plane.** In this section we consider a diffusion process $(\xi(t), \eta(t); t \geq 0)$ taking values in the half plane $\{(x,y) \in \mathbf{R}^2 : y \geq 0\}$ which evolves as a planar Brownian motion away from the boundary and which is slowly reflecting or sticky at $y = 0$. Such a process may be specified by data $(a_0, \theta_0)$ where $a_0 > 0$ and $\theta_0 > 0$ are positive real parameters. The process $(\xi, \eta)$ then satisfies:

(7.1) $\quad \xi(t) \quad \text{and} \quad \eta(t) - \theta_0 \int_0^t \mathbf{1}_{(\eta(s)=0)} \, ds$

are orthogonal martingales relative to the natural filtration of $(\xi, \eta)$, with quadratic variations given by

(7.2) $\quad \langle \xi \rangle(t) = \int_0^t \mathbf{1}_{(\eta(s)>0)} \, ds + a_0 \int_0^t \mathbf{1}_{(\eta(s)=0)} \, ds,$

(7.3) $\quad \langle \eta \rangle(t) = \int_0^t \mathbf{1}_{(\eta(s)>0)} \, ds.$

These properties characterize a strong Markov diffusion in the half-plane, a fact that may be verified by bare hands or by appealing to the general theory of diffusions with boundary conditions. We will write $\mathbf{P}_{(x,y)}$ for the governing probability measure when the process starts from $(\xi(0), \eta(0)) = (x, y)$.

In the case $a_0 = 1$, the processes $(\xi(t); t \geq 0)$ and $(\eta(t); t \geq 0)$ are independent, the former being a Brownian motion on $\mathbf{R}$ and the latter a sticky Brownian motion on the half line $[0, \infty)$ with parameter $\theta_0$. The case of a general value of $a_0$ may be transformed to this by a time-change. Specifically if $(\xi, \eta)$ is specified by data $(a_0, \theta_0)$, then the process $(\tilde{\xi}, \tilde{\eta})$ given by

$$\tilde{\xi}(A(t)) = \xi(t) \quad \text{and} \quad \tilde{\eta}(A(t)) = \eta(t), \tag{7.4}$$

where $A(t) = \langle \xi \rangle(t)$, is of the same class, specified by data $(\tilde{a}_0, \tilde{\theta}_0) = (1, \theta_0/a_0)$. In what follows we will be interested in the exit distribution of the diffusion $(\xi, \eta)$ from certain domains, and since exit distributions are not altered by time changes it will be enough in the proofs to consider the case $a_0 = 1$.

Let $S(\varepsilon)$ denote the infinite strip of width $\varepsilon > 0$,

$$S(\varepsilon) = \{(x, y) \in \mathbf{R}^2 : 0 < x < \varepsilon, y \geq 0\}, \tag{7.5}$$

and let $\tau_{S(\varepsilon)}$ be the first exit time of $(\xi, \eta)$ from this strip. We are interested in the probability that the process exits via the sticky boundary. It is convenient to adopt coordinates $(\phi_1, \phi_2)$ for the starting point $(x, y)$ given by

$$\tan(\phi_1) = \frac{y}{x} \quad \text{and} \quad \tan(\phi_2) = \frac{y}{\varepsilon - x}. \tag{7.6}$$

PROPOSITION 7.1. *For starting points $(x, y) \in S(\varepsilon)$ we have*

$$\frac{2}{\pi} \max(\phi_1, \phi_2) \leq \mathbf{P}_{(x,y)}(\eta(\tau_{S(\varepsilon)}) \neq 0) \leq \frac{2}{\pi}(\phi_1 + \phi_2) + \frac{2}{\sqrt{\pi}} \frac{\theta_0}{a_0} \varepsilon.$$

PROOF. Assume first that $a_0 = 1$, so the processes $\eta$ and $\xi$ are independent. Then for $x \in (0, \varepsilon)$ we have

$$\mathbf{P}_{(x,0)}(\eta(\tau_{S(\varepsilon)}) = 0) = \mathbf{E}_{(x,0)}[f(\tau_{S(\varepsilon)})],$$

where $f(t) = \mathbf{P}_{(x,0)}(\eta_t = 0)$. The transition probabilities of the sticky Brownian motion $\eta$, given in [2], are such that

$$f(t) = \exp(2t\theta_0^2) \operatorname{erfc}(\sqrt{2t}\theta_0),$$

where $\operatorname{erfc}(x) = \frac{2}{\sqrt{\pi}} \int_x^\infty e^{-z^2} dz$. Using standard estimates on erfc we verify that

$$f''(t) = 4\theta_0^4 f(t) - 2\theta_0^3 \sqrt{\frac{2}{\pi t}} + \theta_0 \sqrt{\frac{1}{2\pi t^3}} \geq 0 \quad \text{for all } t \geq 0,$$

and so $t \mapsto f(t)$ is convex. Consequently we may apply Jensen's inequality to obtain

$$\mathbf{P}_{(x,0)}(\eta(\tau_{S(\varepsilon)}) = 0) \geq f(\mathbf{E}_{(x,0)}[\tau_{S(\varepsilon)}]) = f(x(\varepsilon - x)) \geq f(\varepsilon^2/2) \geq 1 - \frac{2}{\sqrt{\pi}}\theta_0\varepsilon.$$

In view of the time-change argument mentioned above, we have for a general $(a_0, \theta_0)$,

$$\mathbf{P}_{(x,0)}(\eta(\tau_{S(\varepsilon)}) \neq 0) \leq \frac{2}{\sqrt{\pi}}\frac{\theta_0}{a_0}\varepsilon.$$

Now, in order to consider a general starting point $(x, y)$, introduce the three stopping times, $\tau_0 = \inf\{t \geq 0 : \eta(t) = 0\}$, $\tau_1 = \inf\{t \geq 0 : \xi(t) = 0\}$ and $\tau_2 = \inf\{t \geq 0 : \xi(t) = \varepsilon\}$. By a standard result on planar Brownian motion

$$\mathbf{P}_{(x,y)}(\tau_1 < \tau_0) = \frac{2}{\pi}\phi_1 \quad \text{and} \quad \mathbf{P}_{(x,y)}(\tau_2 < \tau_0) = \frac{2}{\pi}\phi_2.$$

The lower bound the statement of the proposition follows from this and the relationship between events

$$\{\tau_1 < \tau_0\} \cup \{\tau_2 < \tau_0\} \subseteq \{\eta(\tau_{S(\varepsilon)}) \neq 0\}.$$

For the upper bound we observe that

$$\{\eta(\tau_{S(\varepsilon)}) \neq 0\} \subseteq \{\tau_1 < \tau_0\} \cup \{\tau_2 < \tau_0\} \cup \{\eta(\tau_{S(\varepsilon)}) \neq 0 \text{ and } \tau_0 < \min(\tau_1, \tau_2)\}$$

and use the strong Markov property at time $\tau_0$ together our previous upper bound for $\mathbf{P}_{(x,0)}(\eta(\tau_{S(\varepsilon)}) \neq 0)$ to estimate the probability of the third event on the right-hand side. □

Fix two acute angles $\bar{\phi}_1$ and $\bar{\phi}_2$ and consider the triangle given by

(7.7)
$$\triangle(\varepsilon) = \{(x, y) \in \mathbf{R}^2 : 0 < x < \varepsilon,$$
$$0 \leq y < x\tan(\bar{\phi}_1), 0 \leq y < (\varepsilon - x)\tan(\bar{\phi}_2)\}.$$

Let $\tau_{\triangle(\varepsilon)}$ be the first exit time of $(\xi, \eta)$ from this triangle.

PROPOSITION 7.2. *For starting points $(x, y) \in \triangle(\varepsilon)$ we have*

$$\mathbf{P}_{(x,y)}(\eta(\tau_{\triangle(\varepsilon)}) \neq 0) \leq C(\bar{\phi}_1, \bar{\phi}_2)\left(\phi_1 + \phi_2 + \frac{\theta_0}{a_0}\varepsilon\right),$$

*where $C(\bar{\phi}_1, \bar{\phi}_2)$ depends on $\bar{\phi}_1$ and $\bar{\phi}_2$ alone.*

PROOF. Define $h(x, y) = \mathbf{P}_{(x,y)}(\eta(\tau_{S(\varepsilon)}) \neq 0)$ for $(x, y) \in S(\varepsilon)$. Then, for starting points $(x, y) \in \triangle(\varepsilon)$, applying the strong Markov property at time $\tau_{\triangle(\varepsilon)}$ we obtain

$$h(x, y) = \mathbf{E}_{x,y}[h(\tau_{\triangle(\varepsilon)})\mathbf{1}_{(\eta(\tau_{\triangle(\varepsilon)}) \neq 0)}].$$

Now we estimate $h(x,y)$ using the upper bound from the previous proposition, and $h(\tau_{\triangle(\varepsilon)})$ using the lower bound, to obtain

$$\frac{2}{\pi}(\phi_1 + \phi_2) + \frac{2}{\sqrt{\pi}}\frac{\theta_0}{a_0}\varepsilon \geq \frac{2}{\pi}\min(\bar{\phi}_1, \bar{\phi}_2)\mathbf{P}_{(x,y)}(\eta(\tau_{\triangle(\varepsilon)}) \neq 0)$$

from which the claimed result follows. $\square$

**8. Stochastic flows of kernels.** Suppose that $(K_{s,t}; s \leq t)$ is a stochastic flow of kernels on a space $(E, \mathcal{E})$ as described in the Introduction. A powerful approach to describing $K$ is by means of its family of $N$-point motions. For each integer $N \geq 1$, the $N$-point motion of the flow is a Markov process on $E^N$. Formally it may be described by means of its semigroup which is given by

(8.1) $\quad P_t^N(x, A) = \mathbf{E}[K_{0,t}(x_1, A_1)K_{0,t}(x_2, A_2)\cdots K_{0,t}(x_N, A_N)]$

for all $x = (x_1, x_2, \ldots, x_N) \in E^N$ and $A = A_1 \times A_2 \times \cdots \times A_N \in \mathcal{E}^N$. Informally it should be thought of as describing the motion of $N$ infinitesimal particles sampled from the flow of mass, or if $K$ is interpreted as a random environment governing the evolution of a particle, then take $N$ such particles and let them evolve conditionally independently given $K$. Notice that the family of $N$-point motions is consistent in that any $M$ coordinates (regardless of order) taken from the $N$-dimensional process are distributed as the $M$-dimensional process in the family. The law of the flow $K$ (in the sense of finite-dimensional distributions) is uniquely determined by the associated family of $N$-point semigroups $(P_t^N; t \geq 0)$ for $N \geq 1$. In the opposite direction, in order for a flow of kernels (associated with a given family of $N$-point motions) to exist, we need topological ingredients. Suppose that $E$ is a compact metric space. Then Le Jan and Raimond (Theorem 2.1 of [10]) have proved that whenever $(P_t^N; t \geq 0)$ is a consistent (called compatible there) family of Feller semigroups, there exists an associated flow of kernels; see also Sections 7g and 7h of [16]. This existence result also holds if $E$ is a locally compact separable metric space.

Let $(\theta(k:l); k, l \geq 0)$ be a family of parameters satisfying the hypotheses of Theorem 2.1. Existence and uniqueness of the solution the $\mathcal{A}_N^\theta$-martingale problem gives us a consistent family of semigroups $(P_t^{N,\theta}; t \geq 0)$ for $N \geq 1$. In order to assert the existence of an associated flow $(K_{s,t}^\theta; s \leq t)$ we need the Feller property. This is proved by the same coupling technique as used in [11], Theorem 7.

PROPOSITION 8.1. *For each $N \geq 1$ the semigroup $(P_t^{N,\theta}; t \geq 0)$ associated with $\mathcal{A}_N^\theta$-martingale problem has the Feller property.*

Recall from the Introduction that to each probability measure $\mu$ on $[0,1]$ we associated a stochastic flow of kernels $(K_{s,t}^{\mu}; s \leq t)$ on $\mathbf{Z}$.

PROPOSITION 8.2. *The $N$-point motion of the flow $K^{\mu}$ is a Markov chain on $\mathbf{Z}^N$ with generator $\mathcal{G}_N^p$ given by (6.2) where*

$$p(k\!:\!l) = \int x^k (1-x)^l \mu(dx).$$

PROOF. Fix a realization of the environment $(\Lambda, R)$ and consider $N$ independent particles, each moving as described in the Introduction. Suppose that the configuration of the particles is $x \in \mathbf{Z}^N$ at time $t$, and that the point $(t,y) \in \Lambda$ where $y \in \mathbf{Z}$ is the position of at least one of the particles. Then the configuration jumps to $x + v_{IJ}$ for some $v_{IJ} \in \mathcal{V}(x)$ such that $I \cup J = \{k\!:\!x_k = y\}$ with probability $r^{|I|}(1-r)^{|J|}$ where $r = R(t,y)$. Integrating over possible environments, we find the $N$-point motion jumps from $x$ to $x + v_{IJ}$ for $v_{IJ} \in \mathcal{V}(x)$ at rate $p(|I|,|J|)$. □

Now suppose that $(\mu_n; n \geq 1)$ is a sequence of probability measures on $[0,1]$, with which we associate flows $K^{\mu_n}$ on the integer lattice $\mathbf{Z}$. Assume each $\mu_n$ is centered, in that $\int_0^1 x\mu_n(dx) = 1/2$. Let $\tilde{K}^{\mu_n}$ be the flow of kernels on the scaled lattice $n^{-1/2}\mathbf{Z}$ satisfying

(8.2) $$\tilde{K}_{s,t}^{\mu_n}(x, A) = K_{ns,nt}^{\mu_n}(n^{1/2}x, n^{1/2}A).$$

In view of the convergence results on $N$-point motions established in Section 6, it is natural for us to look for convergence of these scaled flows. But before stating the result we must clarify what is meant by convergence in distribution for stochastic flows of kernels. We will adopt the following definition, similar to that proposed in [9]. Suppose $K^n$, for $n \geq 1$, is a flow of kernels on the scaled lattice $n^{-1/2}\mathbf{Z}$ and $K$ is a flow of kernels on $\mathbf{R}$. We will say $K^n$ converges in distribution to $K$ if the following to holds. Suppose that for each $n \geq 1$ we have a probability measure $\lambda_n$ on the scaled lattice $n^{-1/2}\mathbf{Z}^N$, such that $\lambda_n$ converges weakly to a probability measure $\lambda$ on $\mathbf{R}^N$ as $n$ tends to infinity. Then the $N$-point motion $(Y^n(t); t \geq 0)$ associated with the flow $K^n$, and with the distribution of $Y^n(0)$ being given by $\lambda_n$, must converge in distribution to the $N$-point motion $(Y(t); t \geq 0)$ associated with $K$, with the distribution of $Y(0)$ being given by $\lambda$.

THEOREM 8.1. *Suppose the sequence of centered probability measures $(\mu_n; n \geq 1)$ is such that, as $n$ tends to infinity,*

$$\sqrt{n}x(1-x)\mu_n(dx) \qquad \text{converges weakly to } \nu(dx),$$

*where $\nu$ is some finite measure on $[0,1]$. Then as $n$ tends to infinity the sequence of scaled flows $\tilde{K}^{\mu_n}$ converges in distribution to a flow $K^{\theta}$ on $\mathbf{R}$*

whose $N$-point motions of $K^\theta$ solve the $\mathcal{A}_N^\theta$-martingale problem for the consistent family $(\theta(k\!:\!l); k,l \geq 0)$ determined by

$$\theta(k\!:\!l) = \int_0^1 x^{k-1}(1-x)^{l-1}\nu(dx) \qquad \text{for } k,l \geq 1,$$

and $\theta(1\!:\!0) = -\theta(0\!:\!1) = 0$.

PROOF. Defining $p_n(k\!:\!l) = \int x^k(1-x)^l \mu_n(dx)$, the stated weak convergence to $\nu$, together with the fact that the $\mu_n$ are centred, implies that (6.6) holds, where $\theta$ is as specified in the statement of the theorem. Now fix some $N$, and for $n = 1, 2, \ldots$ let $\lambda_n$ be a probability measure on the scaled lattice $n^{-1/2}\mathbf{Z}^N$. Suppose that as $n$ tends to infinity $\lambda_n$ converges weakly to a probability measure $\lambda$ on $\mathbf{R}^N$. As in Section 6, let $Y^n$ be the scaled process given by $Y^n(t) = n^{-1/2}Y(nt)$ for $t \geq 0$, where $Y$ is a Markov chain with generator $\mathcal{G}_N^{p_n}$, only this time let $Y^n(0)$ have the distribution $\lambda_n$. The arguments of Section 6 still apply, and show that if any subsequence of $Y^n$ converges in distribution to a process $Y$, then $Y$ solves the $\mathcal{A}_N^\theta$-martingale problem starting from $\lambda$. Moreover the uniqueness of the solution of $\mathcal{A}_N^\theta$-martingale problem starting from any point $x \in \mathbf{R}^N$ implies uniqueness starting from $\lambda$ also. The processes $Y^n$ are tight, and thus it follows that they converge in distribution to $Y$, the unique solution to the $\mathcal{A}_N^\theta$-martingale problem starting from $\lambda$.  □

We end the paper by making some comments on the relationship between the flows just constructed, and those of Le Jan and Raimond [10]. The latter flows have $N$-point motions which are symmetric and characterized by means of Dirichlet forms. They are constructed on the circle rather than the real line, but it is reasonable to suppose that there are corresponding flows on the real line with the same local behavior. These flows will belong to the class that we have constructed, and are thus determined by some choice of the measure $\nu$. By straightforwardly applying Itô's formula it is possible to extend the generator $\mathcal{A}_N^\theta$ to functions which are finite linear combinations of products $f_1 f_2$ with $f_1 \in L_N$ and $f_2 \in C^2(\mathbf{R}^N)$. Then, in the case that $\nu$ is a constant multiple of Lebesgue measure on $[0,1]$, by integrating by parts we are able to verify that $\mathcal{A}_N^\theta$ determines a quadratic form having the same structure as those appearing in [10]. This leads us to believe that it is this choice of $\nu$ that gives the Le Jan–Raimond sticky flows.

DEPARTMENT OF STATISTICS
UNIVERSITY OF WARWICK
COVENTRY CV4 7AL
UNITED KINGDOM
E-MAIL: christopherhowitt@hotmail.com
j.warren@warwick.ac.uk